\documentclass{article}

\usepackage[utf8]{inputenc}
\usepackage[T1]{fontenc}
\usepackage{graphicx}
\usepackage{url}
\usepackage{xcolor}
\usepackage[ruled,vlined]{algorithm2e}
\usepackage{setspace}
\usepackage{amsmath}
\usepackage{amssymb}
\usepackage{mathtools} 
\usepackage{siunitx}
\usepackage{makecell} 
\usepackage{todonotes}
\usepackage{xifthen}
\usepackage{hyperref}
\usepackage{array} 
\usepackage{subcaption} 
\usepackage[normalem]{ulem}
\usepackage[american]{babel}
\usepackage{csquotes}
\usepackage[style=apa,dashed=false,doi=false,url=false,annotation=false]{biblatex} 
\usepackage{twoopt}
\usepackage{footnote}
\usepackage{xcolor}
\usepackage{authblk}

\renewcommand{\cite}[1]{\parencite{#1}}

\newcommand{\dfnote}[2][]{%
    {\color{olive}#1}%
    \ifthenelse{\isempty{#2}{}}{}{\todo[color=green,inline]{[DF] #2}}%
}

\definecolor{revcolor}{RGB}{1, 145, 39}
\newcommand{\rev}[1]{#1}

\definecolor{secondrevcolor}{RGB}{3, 0, 235}
\newcommand{\srev}[1]{#1}

\definecolor{thirdrevcolor}{RGB}{232, 93, 0}
\newcommand{\trev}[1]{#1}

\definecolor{fourthrevcolor}{RGB}{171, 52, 235}
\newcommand{\frev}[1]{#1}

\definecolor{fifthrevcolor}{RGB}{2, 194, 172}
\newcommand{\firev}[1]{#1}

\newcounter{mycount}
\newcommand\myprob[1]{%
   \stepcounter{mycount}%
   \vspace{2mm}%
   \noindent \textbf{Problem \themycount} (#1)
}

\newcolumntype{+}{>{\global\let\currentrowstyle\relax}}
\newcolumntype{-}{>{\currentrowstyle}}
\newcommand{\rowstyle}[1]{\gdef\currentrowstyle{#1}%
  #1\ignorespaces
}

\newcommand{\feas}{f}
\newcommandtwoopt{\FG}[2][v][]{%
    \Gamma_{#1}\ifthenelse{\isempty{#2}}{}{^{#2}}%
}

\let\oldnl\nl
\newcommand{\nonl}{\renewcommand{\nl}{\let\nl\oldnl}}

\title{Large-scale Online Ridesharing: The Effect of Assignment Optimality on System Performance\thanks{This is an Accepted Manuscript of an article published by Taylor \& Francis in the Journal of Intelligent Transportation Systems on 4th December 2022, available online: \url{https://www.tandfonline.com/doi/abs/10.1080/15472450.2022.2121651}}}

\author[1]{David Fiedler\thanks{CONTACT David Fiedler. Email: david.fiedler@agents.fel.cvut.cz}}
\author[1]{Michal \v{C}ertick\'{y}}
\author[2]{Javier Alonso-Mora}
\author[1]{Michal P\v{e}chou\v{c}ek}
\author[1]{Michal \v{C}\'{a}p}
\affil[1]{Department of Computer Science, Faculty of Electrical Engineering, CTU in Prague, Czech Republic}
\affil[2]{Department of Cognitive Robotics, 3ME, TU Delft, Netherlands}

\begin{document}

\maketitle

\begin{abstract}
Mobility-on-demand (MoD) systems consist of a fleet of shared vehicles that can be hailed for one-way point-to-point trips. The total distance driven by the vehicles and the fleet size can be reduced by employing ridesharing, i.e., by assigning multiple passengers to one vehicle.
However, finding the optimal passenger-vehicle assignment in an MoD system is a hard combinatorial problem.
In this work, we demonstrate how the VGA method, a recently proposed systematic method for ridesharing, can be used to compute the optimal passenger-vehicle assignments and corresponding vehicle routes in a massive-scale MoD system.
In contrast to existing works, we solve all passenger-vehicle assignment problems to {\em optimality}, regularly dealing with instances containing thousands of vehicles and passengers. 
Moreover, to examine the impact of using optimal ridesharing assignments, we compare the performance of an MoD system that uses optimal assignments against an MoD system that uses assignments computed using insertion heuristic and against an MoD system that uses no ridesharing. 
We found that the system that uses optimal ridesharing assignments subject to the maximum travel delay of 4 minutes reduces the vehicle distance driven by \SI{57}{\percent} compared to an MoD system without ridesharing.
Furthermore, we found that the optimal assignments result in a \SI{20}{\percent} reduction in vehicle distance driven and \SI{5}{\percent} lower average passenger travel delay compared to a system that uses insertion heuristic.
\end{abstract}


\section{Introduction}
In densely populated cities, private cars are considered as an unsustainable mode of transportation. 
Typically, parking capacity and road capacity are insufficient to accommodate all private transport and, at the same time, difficult to expand due to lack of available urban space or high cost. 
As a result, many modern cities suffer from traffic congestion, unavailability of parking spaces, and air pollution.

\trev{
One of the proposed solutions to address these problems is the deployment of metropolitan mobility-on-demand (MoD) systems providing an alternative to traveling in a private vehicle designed to be as comfortable as traveling in a private car but with smaller parking capacity and road capacity requirements~\cite{spieserSystematicApproachDesign2014b,millerPredictivePositioningQuality2017,alonso-moraOndemandHighcapacityRidesharing2017,capMultiObjectiveAnalysisRidesharing2018}. 
These MoD systems consist of a fleet of shared passenger vehicles that jointly serve the travel requests of the system's users. 
Usually, these systems are assumed to use traditional vehicles able to carry up to four passengers simultaneously. 
For each incoming travel request, the MoD system assigns the request to one of the vehicles and alters its route such that the passenger is picked up and transported to the drop-off location.
Examples of such an MoD system include services like Uber Pool or Lyft Line, as well as the future systems of autonomous self-driving cars being developed by companies such as Waymo,  Cruise, or Motional. }

Such MoD systems employ {\em vehicle sharing}, so they can serve the existing transportation demand with a smaller, highly-utilized vehicle fleet and thus, significantly reduce the need for urban parking space.
To further improve the system's efficiency, the provider can implement {\em ridesharing}, where multiple passengers can be transported in one vehicle simultaneously\srev{~\cite{alonso-moraOndemandHighcapacityRidesharing2017}}. 
Efficient ridesharing increases vehicle occupancy, which consequently reduces the required fleet size and total distance driven by the vehicle fleet, resulting in ecological and economic benefits. 

\trev{
Two clarifying notes on terminology are in order. 
First, we note that MoD systems are ultimately envisioned also to include high-capacity transportation modes (e.g., trains, subway, buses) and to allow for transfers between different vehicles~\cite{susanshaheenSimilaritiesDifferencesMobility2020}. 
However, this article focuses on MoD systems that transport each passenger from their pick-up to their destination in one vehicle. Analysis of MoD systems that allow transfers is left for future work.
Second, we also caution that the term ridesharing is overloaded. 
This paper focuses on ridesharing in on-demand mobility systems, where each vehicle is driven by a for-hire driver who transports travelers between their desired pick-up and drop-off locations. 
Alternatively, in the future, these vehicles may be self-driving. 
Apart from that, there is a distinct concept called \emph{peer-to-peer ridesharing}, where the vehicle is typically owned and driven by one of the travelers, whose primary motivation is to reach his/her intended destination.
Readers interested in peer-to-peer ridesharing are referred to the growing body of research devoted to this model, for example, the work of ~\textcite{liRestrictedPathbasedRidesharing2019} that studies the impact of high occupancy toll lane configurations on the willingness of (peer) drivers to share a ride, or~\textcite{maRidesharingUserEquilibrium2020} and \textcite{yanStochasticRidesharingUser2019} who study the ridesharing user equilibrium in the context of peer-to-peer ridesharing.
}
%

\subsection{Related work}
\label{sec:related}
Recently, a number of mobility-on-demand system models have been developed with the aim to provide quantitative insights into the potential of large-scale carsharing and ridesharing to improve the efficiency of urban transportation.

 Most existing models of MoD systems assume unit-capacity vehicles~\cite{spieserSystematicApproachDesign2014a,bischoffSimulationCitywideReplacement2016,fiedlerImpactMobilityondemandTraffic2017,maciejewskiCongestionEffectsAutonomous2018,venkatramanCongestionawareTabuSearch2019}. 
However, transportation systems that do not employ ridesharing suffer from poor operational efficiency because the vehicles need to travel empty from the drop-off point of a passenger to the pick-up point of the following passenger.
Such unallocated trips can generate significant extra vehicular traffic in the system; various studies indicate the growth in vehicle distance traveled from \SIrange{17}{40}{\percent} depending on the system configuration~\cite{bischoffSimulationCitywideReplacement2016,fiedlerImpactMobilityondemandTraffic2017,maciejewskiCongestionEffectsAutonomous2018}.
The average vehicle occupancy observed in such systems is considerably lower than one passenger per vehicle~\cite{fiedlerImpactRidesharingMobilityonDemand2018}, a finding which also corresponds to the average vehicle occupancy measured in already operating taxi services~\cite{nyctaxilimousinecommission2016TLCFactbook2016}.
The low occupancy in MoD systems can lead to congestion, which could be partially alleviated by a congestion-aware dispatching~\cite{venkatramanCongestionawareTabuSearch2019}.

Therefore, it is beneficial to consider vehicles with a capacity higher than one and allow ridesharing between passengers. 
%
In contrast to peer-to-peer ridesharing~\cite{masoudRealtimeAlgorithmSolve2017,liRestrictedPathbasedRidesharing2019,tamannaeiCarpoolingProblemNew2019}, here we are interested in the centralized setting, where a central dispatcher decides on an efficient assignment of travel requests to fleet vehicles. This problem is commonly formulated as a Vehicle Routing Problem with Pickup and Deliveries~(VRPPD) or, more specifically, as Dial-a-Ride Problem~(DARP)~\cite{cordeauDialarideProblemModels2007,tothVehicleRoutingProblems2014}.
These formulations can be solved optimally using off-the-shelf Integer Linear Programming (ILP) solvers or domain-tailored ILP solution techniques.
\trev{
However, the applicability of these methods is limited to small-scale instances with at most tens of requests and vehicles. Large-scale MoD systems typically require the ability to find routes for many more vehicles and requests. For example, in New York City (NYC), there are almost \num{100000} active taxis per hour during peak traffic~\cite{nyctaxilimousinecommission2018Factbook2018}.
Therefore, DARP instances appearing in large-scale MoD systems are typically solved using heuristic methods.
}

A popular heuristic method for large-scale dynamic DARP is the Insertion Heuristic (IH)~\cite{campbellEfficientInsertionHeuristics2004,kalinaAgentsVehicleRouting2015,bischoffCitywideSharedTaxis2017,fiedlerImpactRidesharingMobilityonDemand2018}.
Also, IH is often used as a subcomponent of more sophisticated algorithms. For example, in ridesharing with demand prediction~\cite{vanengelenEnhancingFlexibleTransport2018}, when integrating ridesharing with public transport~\cite{maDynamicRidesharingDispatch2019}, or as an initial solution generator for metaheuristic methods~\cite{muelasVariableNeighborhoodSearch2013}.
The \rev{metaheuristic} methods, which are effective in solving conventional DARP problem instances~\cite{hoSurveyDialarideProblems2018}, 
\firev{typically target scenarios with less than twenty vehicles~\cite{masmoudiThreeEffectiveMetaheuristics2016,pfeifferALNSAlgorithmStatic2022} and suffer from scalability issues when applied to large-scale DARPs. 
}
\rev{
However, in the last decade, there was some progress with metaheuristic approaches enabling solving larger instances.
\textcite{jungDynamicSharedTaxiDispatch2015} used simulated annealing to solve scenarios with \num{600} operating vehicles. 
Another popular metaheuristic is the Greedy Randomized Adaptive Search Procedure (GRASP) used by~\textcite{santosDynamicTaxiRidesharing2013}.
The authors were able to solve instances with up to \num{750} requests.
They also tested an online setting with \num{78000} requests per day, and later, they improved the results significantly~\cite{santosTaxiRideSharing2015}. 
}
\textcite{muelasVariableNeighborhoodSearch2013}~also solved four types of specialized DARP scenarios with up to 90 vehicles using Variable Neighborhood Search.
Later, \textcite{muelasDistributedVNSAlgorithm2015}~modified this approach to a distributed version which was able to solve scenarios with up to \num{1668} vehicles and \num{16000} requests.
\firev{Another metaheuristic, a modified artificial bee colony algorithm, was used by~\textcite{zhanModifiedArtificialBee2021}.
The method was able to solve an instance of \SI{3661}{requests} and \SI{2400}{vehicles}.
Later, this method was used in a simulation-optimization framework for an MoD system with electric vehicles~\cite{zhanSimulationOptimizationFramework2022}.
}

A systematic and scalable approach for pairwise ridesharing based on bipartite matching in the so-called shareability network was proposed by~\textcite{santiQuantifyingBenefitsVehicle2014}. 
The analysis revealed that up to \SI{80}{\percent} of the trips could be pairwise shared while keeping the travel delay lower than a couple of minutes. 
Later, \textcite{alonso-moraOndemandHighcapacityRidesharing2017} proposed a new method that lifted the limit of two passengers per car and evaluated this method on the NYC taxi dataset. 
Finally, \textcite{capMultiObjectiveAnalysisRidesharing2018}~utilized this method to study the trade-offs between the quality of service and the operation cost inherent in ridesharing. 

\firev{
Finally, apart from optimizing the assignment of passengers to vehicles, we can also optimize the pickup and drop-off locations if we enable short walking for passengers. 
Such an approach was tested by~\textcite{fielbaumOndemandRidesharingOptimized2021}, showing that it can improve the level of service and decrease the total travel time.
Later, \textcite{fielbaumOptimizingVehicleRoute2021} tried optimizing pickup and drop-off positions of precomputed vehicle plans to measure the benefits exactly.
He demonstrated that we could decrease the travel cost by almost \SI{19}{\percent} when optimizing the locations with a heuristic method and \SI{24}{\percent} with a slower optimal solution method.
}

\subsection{Contribution}
In this work, we extend the existing study of ridesharing in large-scale MoD systems~\cite{fiedlerImpactRidesharingMobilityonDemand2018} by analyzing the impact of passenger-vehicle assignment optimality on system performance. 
To do this, we use a variant of the vehicle-group assignment (VGA) method used by~\textcite{alonso-moraOndemandHighcapacityRidesharing2017} and \textcite{capMultiObjectiveAnalysisRidesharing2018}.  
\trev{We chose this method because it was previously demonstrated to be able to efficiently solve large-scale DARP instances with tight pick-up and drop-off time windows that are characteristic of large-scale MoD systems.
Moreover, it is less complex than the classical exact methods for DARP, and it can be easily modified to a resource-constrained version, which we also evaluate in this work.
}

The contribution of this paper is 3-fold:

\emph{1) Optimality}: We took special care to ensure that all ridesharing assignments and routes are computed optimally. 
This is in contrast to \textcite{alonso-moraOndemandHighcapacityRidesharing2017}, who used a similar solution algorithm to evaluate shareability within the NYC taxi dataset, but to maintain computational tractability, 
\rev{ the actual implementation used in the experiment resorted to heuristics and time-outs, leading to suboptimal performance of the system. Moreover, it remained unclear how far are the reported performance metrics from optimum.}
In this work, we identified and solved several algorithmic bottlenecks, and consequently, we were able to obtain optimal \rev{ridesharing assignments} \trev{for the majority of the evaluated scenarios}. 

\emph{2) Scale}: We implemented performance optimizations that enable us to significantly scale the algorithm and compute optimal ridesharing assignments for instances of unprecedented size peaking at more than \num{21000} active travel requests and \num{11000} vehicles. This is in contrast to~\textcite{capMultiObjectiveAnalysisRidesharing2018} who proposed the optimal version of the VGA method but were only able to solve problem instances with a bit less than 500 requests.

\emph{3) Impact of Assignment Optimality}: 
\rev{
With an 1) {\em optimal} and 2) {\em scalable} implementation of a ridesharing algorithm, we are able to achieve the main objective of this work: to quantify the impact of using an optimal ridesharing method on system performance in comparison to the performance achieved by sub-optimal ridesharing methods.
We quantify the reduction in vehicle distance traveled, travel delay, used vehicles, and traffic density for different ridesharing strategies.
}
Specifically, we compare the above metrics in \rev{six} scenarios: a) a present-day transportation using private vehicles, b) an MoD system without ridesharing, c) an MoD system with ridesharing based on the IH, d) an MoD system with optimal ridesharing computed using Vehicle Group Assignment (VGA) method, e)~an~MoD system with ridesharing solved by a VGA method \rev{with limited computational resources, and f) an~MoD system with ridesharing solved by a resource-constrained variant of VGA method as implemented in~\textcite{alonso-moraOndemandHighcapacityRidesharing2017}}. 
This allows us to give a quantitative answer to the question of how much do we gain by actually taking the effort to compute optimal assignments? 

The performance comparison of the system that uses optimal assignments against the system that uses IH is particularly interesting, as the latter approach is widely used in existing studies of large-scale MoD systems~\cite{campbellEfficientInsertionHeuristics2004,bischoffCitywideSharedTaxis2017,fiedlerImpactRidesharingMobilityonDemand2018}, while the former represents the fundamental bound on system performance.

Our evaluation revealed that optimal ridesharing assignments can reduce the distance driven in the system by \SI{57}{\percent} compared to an MoD system without ridesharing, and simultaneously, we managed to maintain the passenger travel delay below 4 minutes.
Furthermore, we found that the optimal ridesharing assignments are considerably more efficient than the assignments computed by IH. Specifically, in the system that uses optimal assignments, the total vehicle distance driven is reduced by \SI{20}{\percent}, and simultaneously, average passenger travel delay is reduced by~\SI{5}{\percent}.
\frev{
Moreover, in order to provide insights into the limits of the VGA method, we performed a sensitivity analysis on our city-scale scenarios with respect to 1) the length of a ridesharing batch, 2) the vehicle capacity, and 3) and the maximum allowed passenger delay. 
Our results show that the VGA method is capable of finding optimal assignments given that vehicle capacity is no more than 5-10 passengers and the maximum allowed delay is no more than 4-7 minutes, depending on the demand structure and intensity. For scenarios using higher-capacity vehicles or with more permissive delay constraints, the VGA algorithm can no longer certify optimality of the computed ridesharing assignments. 


}


\section{Methodology}
\label{sec:methodology}

We use a travel demand model to generate a dataset of all private car trips in Prague. 
Then, we design an MoD system that can serve these existing trips with the required service quality \rev{(measured by maximum travel delay)}.
After that, we implement the considered solution methods for passenger-vehicle matching. 
Finally, we simulate various scenarios in multi-agent simulation and analyze the results.

\subsection{Input data}
\label{sec:demand}
The set of trips that represent the transportation demand is generated by the multi-agent activity-based model \rev{described in~\textcite{drchalDatadrivenActivityScheduler2019}}.
We chose the city of Prague, the Czech Republic for a case study because a) we have access to the travel demand model for the area and b) because its demand density, demand structure, and road topology are representative for a large European city. 
This is in contrast to previously considered urban areas, such as Manhattan or Singapore, which due to an extremely high density of travel demand, lead to overly-optimistic estimates \rev{of the benefits of ridesharing}.
\rev{Nevertheless, for interested readers, we also performed a version of our experiment using the Manhattan taxi demand dataset (the dataset previously used by \textcite{alonso-moraOndemandHighcapacityRidesharing2017}), and we present the results in Appendix~\ref{sec:manhattan}.}

In contrast to traditional four-step demand models~\cite{hensher2007handbook}, which use trips as the fundamental modeling unit, activity-based models employ so-called activities (e.g., work, shop, sleep) and their sequences to represent the transport-related behavior of the population. 
Travel demand then occurs due to the agents' necessity to satisfy their needs through activities performed at different places at different times. 
These activities are arranged in time and space into sequential daily schedules. Trip origins, destinations, and times are endogenous outcomes of activity scheduling. 
The activity-based approach considers individual trips in context and therefore allows representing realistic trip chains.

The model used in this work covers a typical workday in the metropolitan area of Prague. 
The population of over \SI{1.3}{million} is modeled by the same number of autonomous, self-interested agents, whose behavior is influenced by their sociodemographic attributes, current needs, and situational context. 
\rev{Individual decisions of the agents are implemented using four modules responsible for choosing the activity type, duration, location, and mode.
Each module uses a dedicated machine learning model (such as neural network, decision trees, regression tree) trained so that its output matches various real-world data sets such as travel diaries and other transportation-related surveys, demographic data, points of interests, and transport network structure.}
Planned activity schedules are simulated and tuned, and finally, their temporal, spatial, and structural properties are validated against additional historical real-world data (origin-destination matrices and surveys) using the six-step validation framework VALFRAM~\cite{jass2016, drchal2015data}.

\begin{figure}[ht]
\centering{}\includegraphics[width=1\columnwidth]{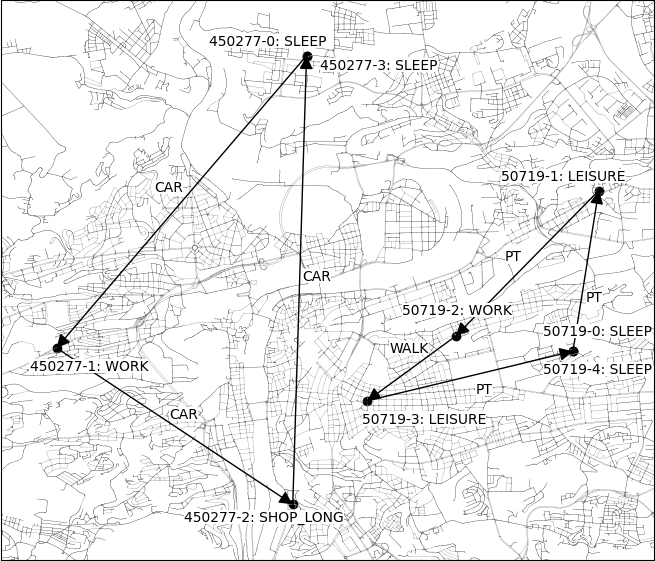}\caption{\label{fig:trips_sample}\rev{Two example trips from the generated demand.
The filled circles represent activities, while the arrows represent trips between those activities. 
Next to each activity, we can see the person and activity IDs in the format \texttt{person\_id-activity\_id}.
The activities corresponding to these IDs can be found in Table~\ref{fig:activities}.}
}
\end{figure}

\begin{figure}[ht]
\centering{}\includegraphics[width=1\columnwidth]{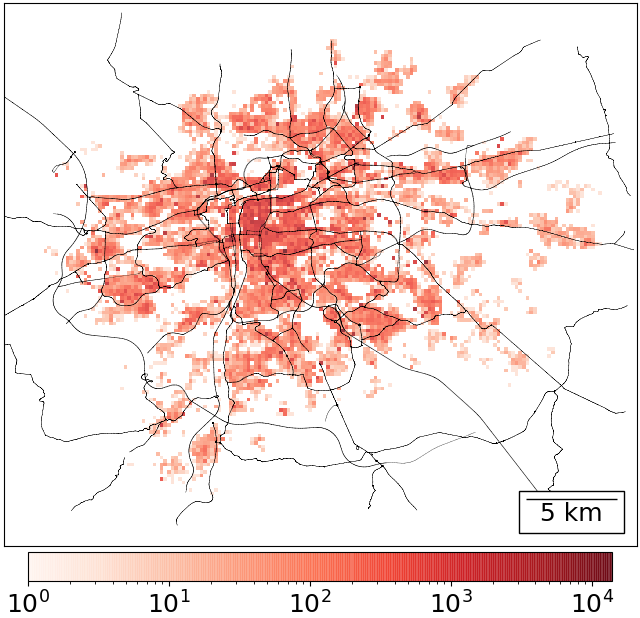}\caption{\label{fig:demand}Demand for personal vehicle traffic in Prague. The start positions of all vehicle trips are discretized to squares of \num{200} square meters. Darker color translates to higher demand, and the color bar has a logarithmic scale.
}
\end{figure}

\begin{figure}[ht]
\centering{}\includegraphics[width=1\columnwidth]{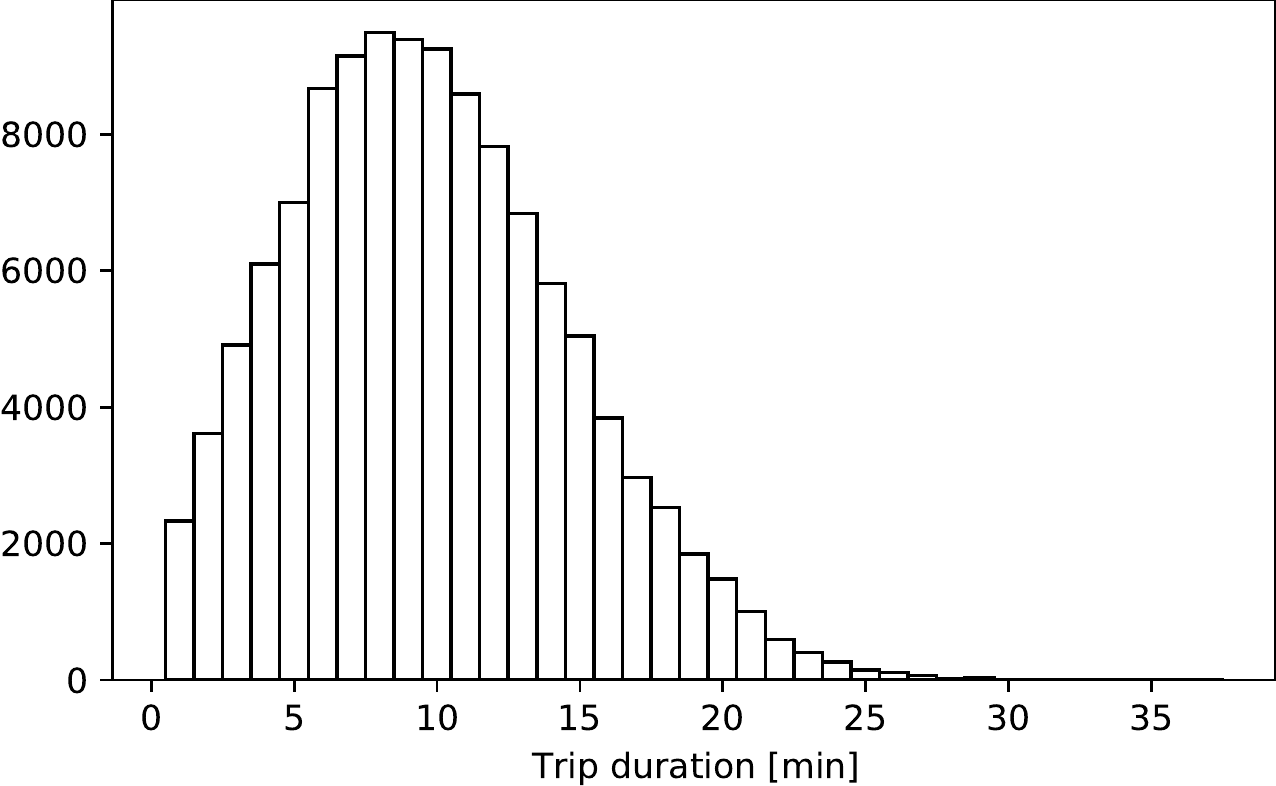}\caption{\label{fig:demand_trip_lengths}\rev{Histogram of fastest path travel times for each trip.}
}
\end{figure}

The model consists of over three million trips by all modes of transport in one 24-hour scenario, out of which there are roughly one million trips realized by private vehicles~(Figure~\ref{fig:demand}).
\rev{Tables~\ref{fig:trips},~\ref{fig:activities},~and~Figure~\ref{fig:trips_sample} show example activity schedules for two agents.}
In this work, we select only the trips realized by private vehicles in two representative time intervals: the {\em peak} dataset includes trips that starts 06:30 and 08:00, and the \emph{off-peak} dataset includes trips that start between 10:30 and 12:00. 
The two datasets contain about \num{130000} and \num{45000} trips, respectively.
\rev{The duration of trips ranges from 1 to 37 minutes; the histogram is in Figure~\ref{fig:demand_trip_lengths}.}

\begin{table}
\scriptsize
\centering{}%
\setlength{\tabcolsep}{0.3em}
{\renewcommand{\arraystretch}{1.2}%
\begin{tabular}{|r|r|r|r|r|}
\hline
\thead{Trip} & \thead{Person} & \thead{From} & \thead{To} & \thead{Mode}
\tabularnewline
\hline
\hline
\num{4500942} & \num{50719} & \num{0} & \num{1} & PT 
\tabularnewline
\hline
\num{4500943} & \num{50719} & \num{1} & \num{2} & PT 
\tabularnewline
\hline
\num{4500944} & \num{50719} & \num{2} & \num{3} & WALK 
\tabularnewline
\hline
\num{4500945} & \num{50719} & \num{3} & \num{4} & PT 
\tabularnewline
\hline
\num{4789903} & \num{450277} & \num{0} & \num{1} & CAR 
\tabularnewline
\hline
\num{4789904} & \num{450277} & \num{1} & \num{2} & CAR 
\tabularnewline
\hline
\num{4789905} & \num{450277} & \num{2} & \num{3} & CAR 
\tabularnewline
\hline
\end{tabular}}

\caption{\label{fig:trips} \rev{Example trips. Each trip connects two activities, shown in Table~\ref{fig:activities}. We can identify each activity by \texttt{Person} column and \texttt{From}/\texttt{To} columns that correspond to the \texttt{Activity} column in Table~\ref{fig:activities}}}
\end{table}

\begin{table}
\scriptsize
\centering{}%
\setlength{\tabcolsep}{0.3em}
{\renewcommand{\arraystretch}{1.2}%
\begin{tabular}{|r|r|r|r|r|r|r|}
\hline
\thead{Person} & \thead{Activity} & \thead{Start} & \thead{End} & \thead{Type} & \thead{Lat} & \thead{Lon}
\tabularnewline
\hline
\hline
\num{50719} & \num{0} & 00:00 & 04:06 & SLEEP & \num{50.084294} & \num{14.490635}
\tabularnewline
\hline
\num{50719} & \num{1} & 06:56 & 07:42 & LEISURE & \num{50.110286} & \num{14.496852}
\tabularnewline
\hline
\num{50719} & \num{2} & 10:31 & 13:49 & WORK & \num{50.086623} & \num{14.461201}
\tabularnewline
\hline
\num{50719} & \num{3} & 15:01 & 16:04 & LEISURE & \num{50.076027} & \num{14.439032}
\tabularnewline
\hline
\num{50719} & \num{4} & 17:17 & 00:00 & SLEEP & \num{50.084294} & \num{14.490635}
\tabularnewline
\hline
\num{450277} & \num{0} & 00:00 & 07:22 & SLEEP & \num{50.131751} & \num{14.423139}
\tabularnewline
\hline
\num{450277} & \num{1} & 07:54 & 15:43 & WORK & \num{50.084170} & \num{14.360924}
\tabularnewline
\hline
\num{450277} & \num{2} & 16:00 & 16:35 & SHOP\_LONG & \num{50.059205} & \num{14.420547}
\tabularnewline
\hline
\num{450277} & \num{3} & 17:26 & 00:00 & SLEEP & \num{50.131751} & \num{14.423139}
\tabularnewline
\hline
\end{tabular}}

\caption{\label{fig:activities} \rev{Example activities.}}
\end{table}

\subsection{System Model}
\label{sec:system_model}

\rev{For MoD systems design, we adopt a station-based methodology described by~\textcite{pavoneRoboticLoadBalancing2012} or by \textcite{wallarOptimizingVehicleDistributions2019a}, which means that idle MoD vehicles are parked in dedicated parking facilities instead of parking on-street or cruising.
This setup is typical in carsharing or bike-sharing systems because curb parking would take valuable urban space, and cruising for parking would increase fuel consumption and congestion.
Further, in case of electric vehicles, stations will provide charging infrastructure.}
Vehicles are initialized in stations and leave a station only to serve travel requests. 
Whenever a vehicle becomes idle, it starts driving to the nearest station to park there. 

\subsubsection{\rev{Station Positioning, Rebalancing, and Fleet-sizing}}
\label{sec:station_pos_reb_and_fleet_sizing}
\begin{figure}
\centering{}\includegraphics[width=1\columnwidth]{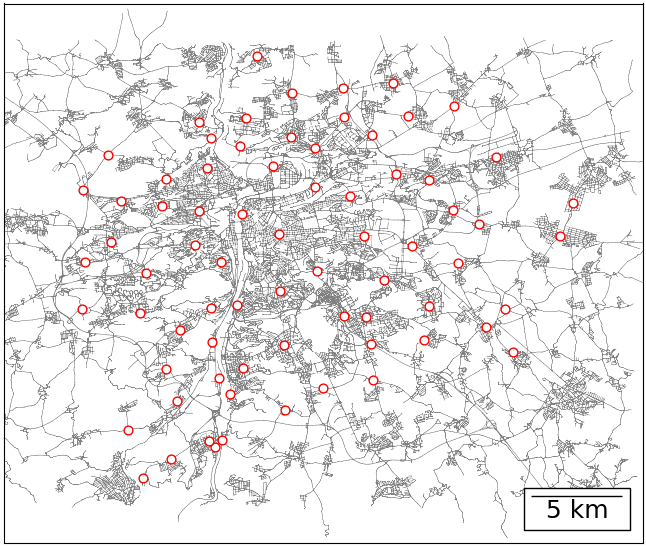}\caption{\label{fig:stations}MoD system stations in the city of Prague. There are 73 stations in total, shown as red circles.
}
\end{figure}

We use 73 stations shown in Figure~\ref{fig:stations} chosen such that every node on the road network (excluding roads without travel requests such as tunnels or highways) can be reached from one of the stations within \rev{\SI{210}{\s}, and the number of stations is minimized. 
We compute the station positions using an integer program with binary variables $ s_n $ for each node $ n $ in the set of serviced nodes $ N $, where each variable $ s_n $ indicates if there is a station at node $ n $ (1) or not (0).
We minimize
\begin{equation}
    \sum_{n \in N} s_n, 
\end{equation}
subject to
\begin{equation}
    \sum_{n' \in P_n} s_{n'}  \geq 1 \quad \forall n \in N,
\end{equation}
where $ P_n $ is a set of nodes from which $ n $ is reachable within \SI{210}{\s}.
}

The stock of vehicles at each station is stabilized by a vehicle rebalancing process that continuously sends empty vehicles from stations with a surplus of vehicles to stations that have a shortage of vehicles. 
We use the rebalancing policy introduced by~\textcite{pavone_robotic_2012} and later evaluated by~\textcite{spieserSystematicApproachDesign2014a} in the Singapore MoD case study.
\rev{In one-minute intervals, we generate an integer program for transferring vehicles from stations with more vehicles compared to the initial state to stations with fewer vehicles compared to the initial state such that the number of vehicles in each station $ s $ is kept above a corresponding threshold $ \tau_s $, and the total length of all rebalancing trips is minimized.
We experimentally determined that in order to compensate for driving vehicles, the $ \tau_s $ should be no more than \SI{85}{\percent} of the initial number of vehicles parked in $ s $.
Also, we use only stations with at least \SI{5}{\percent} more vehicles over the initial state as source stations in order to prevent rebalancing instabilities, i.e., rebalancing flows in the opposite directions.
}

Our objective is to achieve full service availability during the entire experiment, i.e., every request should be served.
\rev{
We experimentally determined that in order to be able to serve every request during the morning peak (see Section~\ref{sec:demand}) without ridesharing, the MoD system requires a total of \num{68201} vehicles\footnote{\rev{The number of vehicles is smaller than the number of trips. This is possible because even without ridesharing, one vehicle can serve more travel requests sequentially.}}.
We used the same fleet size for other scenarios (off-peak, ridesharing)
\footnote{\rev{In practice, we can use lot fewer vehicles, especially for the ridesharing scenarios.
However, since the fleet-sizing problem is not the focus of this article, we used this fleet-sizing method to ensure that the size of the fleet is not the limiting factor. 
The number of vehicles that were actually used in each experiment is in our experimental results.
}} 
as these experiments are guaranteed to require fewer or equal vehicles than the scenario without ridesharing.}
Specifically, to determine the number of vehicles in each station needed to ensure full service availability, we first created a dedicated vehicle for each request in the station closest to the requested pickup location. 
Then, we started iteratively reducing the number of vehicles \rev{by the same factor} in each station until the first vehicle shortage event occurred \rev{in any station}.
\rev{Then, we used the vehicle counts from the last iteration without any shortage.}
This procedure guarantees that there is a sufficient number of vehicles to serve all requests from the nearest station.

\subsubsection{Problem Formulation}
The set of all vehicles will be denoted as $ V=1,\ldots,m $\rev{, all vehicles have the same capacity $ K $}. 
Travel requests are modelled as a sequence $(t_1, o_1, d_1), (t_2, o_2, d_2), \ldots \ $, where $t_i, o_i,$ and $d_i$ are the announcement time, origin point, and destination point of request $i$, respectively. The $i$-th request is revealed only at time $t_i$.
\rev{We obtain requests from the demand model trips simply by setting $ t $ equal to the trip start time, $ o $ equal to the trip start location, and $ d $ equal to the trip end location.}

The state of a vehicle $ v $ at a particular time point encodes its current position, the set of passengers currently on-board of the vehicle, and its current plan.
The plan of a vehicle is represented as a sequence of locations $ p = l_1, l_2,\ldots \ $, where each location $ l_{i} $ is either an origin location~$ o_i $, or a destination location~$ d_i $ of request~$ i $ that is scheduled to be serviced by the plan.
A vehicle plan is \emph{valid} only if the plan contains origin location and later destination location for each onboard passenger.

The operational cost of vehicle $ v $ when following plan $ p $ is denoted $ c( p, v ) $.
For simplicity, we define $ c( p, v )$ to be equal to the distance driven by the vehicle when it follows plan~$ p $.
\rev{Each plan $ p $ requires a vehicle of capacity $ K(p) $, where $ K (p) \leq |p| $.}
The travel delay of request $ r $ when it is served by vehicle $v$ following plan $ p $ is computed as:
\begin{equation}
q_r(p, v) \coloneqq (t_r^\mathrm{dropoff} - t_r) - \delta^\mathrm{baseline}_r.
\end{equation}
Here, $ t_r^\mathrm{dropoff} $ is the time when the request is dropped off under plan $ p $ and $ \delta^\mathrm{baseline}_r $ is the duration along direct route from the request's origin to its destination.
\rev{Note that the passenger's waiting time is included in the delay, and therefore, the maximum delay also limits the maximum waiting time.}

Our goal is to minimize the total operational cost of the system, such that the \rev{delay} of every passenger is bounded by a constant  $ q_\mathrm{max}, $ \rev{and the maximum capacity $ K $ is respected for all vehicles}.
That is, we desire to minimize

\begin{equation}
    \sum_{v \in V} c(p, v)
\end{equation}

subject to

\rev{
\begin{alignat}{2}
    q_r(p, v) &\leq q_\mathrm{max} \quad &&\forall r \in R \\
    K(p_v) &\leq K \quad &&\forall v \in V.
\end{alignat}
}

\subsection{Request-vehicle Matching}
\label{sec:matching}
In an MoD system, new requests dynamically arrive and need to be served.
A ridesharing algorithm tries to find the \emph{optimal system plan} (i.e., a collection of vehicle plans), 
such that 1) every request is served, 2) maximum discomfort constraint $ q_\mathrm{max} $ is respected, and 3) the total operation cost is minimized.
This planning procedure is repeated periodically, and each such planning period is referred to as a \emph{batch}. 
During one batch, we collect all newly announced requests and execute a planning procedure that computes request-vehicle matching and corresponding vehicle plans.
\rev{
We make the following assumptions:~a) travel time on each road segment is constant over time and does not depend on the number of vehicles on the segment, b) the execution of the vehicle schedule is perfect (there are no random delays), and~c) the mode choice is fixed in the demand model and customers accepts any plan that satisfies the max delay constraint (which is guaranteed in our setup, as explained in Section~\ref{sec:station_pos_reb_and_fleet_sizing}). 
}

The request-vehicle matching can be modeled as a Dial-a-Ride (DARP) problem, which is known to be NP-hard~\cite{tothVehicleRoutingProblems2014}.
In this work, we implement and compare two methods for computing such request-vehicle matching.
First, we implement Insertion Heuristic (IH)~\cite{campbellEfficientInsertionHeuristics2004}, a popular heuristic algorithm for DARP and other vehicle routing problems.
Second, we implement Vehicle-Group Assignment (VGA) method, \cite{capMultiObjectiveAnalysisRidesharing2018}, which is a recently proposed exact solution method for DARP exhibiting good scalability properties. 

\subsubsection{Insertion Heuristic}
The pseudocode of the IH is presented in Algorithm~\ref{alg:insertion_heuristic}.
The algorithm is implemented as follows:
For each new request, the IH algorithm attempts to insert the request into the plan of \rev{each} vehicle.
The current plan of a vehicle  $ v $, denoted as $ p_v $, is the plan computed in one of the previous iterations of the algorithm.
For a particular vehicle $ v $, we try all possible indexes $ i $ in plan $ p_v $ to insert pickup of the new request before $ i $ and all possible indexes $ j, j > i $ to insert drop off of the new request before $ j $.
We denote such plan as $ p_v^{\mathrm{new}} $.
Note that the relative ordering of all locations from $ p_v $ remains unchanged in the new plan, \rev{and therefore, optimality is not guaranteed}.
Finally, among all plans generated this way, we select the plan (and the corresponding vehicle) that minimizes the increase in operating cost and at the same time satisfies the service discomfort constraints.

\SetKwProg{alg}{Algorithm}{}{}
\SetKwProg{on}{On}{}{}
\SetKwProg{function}{Function}{}{}
\SetKw{And}{and}
\SetKwInOut{Input}{input}
\SetKwInOut{Output}{output}

\SetKwFunction{pp}{PP}
\LinesNumbered

\begin{algorithm}[t]

\Input{Current plan $ p_v $ of each vehicle $ v $ that was computed in one of the previous iterations of the algorithm and the set of new requests $ D_n $, i.e., requests announced in the $ n $th batch.}
\For{$ r \in D_n $} {
	$ \delta_c^{min} \leftarrow \infty $ \tcc*{min. cost increment}
	$ v^* \leftarrow $ null \;
	\For{$ v \in V $}{
		\For{$ i \in {1, \ldots, |p_v|} $}{
			\For{$ j \in {\rev{i} + 1, \ldots, |p_v| + 1} $}{
		   		$ p_v^{\mathrm{new}} \leftarrow p_v $\;
		   		insert $ o_r $ to $ p_v^{\mathrm{new}} $ before index $ i $\;
		   		insert $ d_r $ to $ p_v^{\mathrm{new}} $ before index $ j $\;
		   		$ \delta_c \leftarrow c(p_v^{\mathrm{new}}) - c(p_v) $\;
		   		\If{$ p_v^{\mathrm{new}} $ is feasible \And $ \delta_c < \delta_c^{min} $ }{
					$ \delta_c^{min} \leftarrow \delta_c $\;
					$ p^* \leftarrow p_v^{\mathrm{new}} $\;
					$ v^* \leftarrow v $\;
				}
			}
		}
	}
	\If{$ v^* $ not null }{
		vehicle $ v^* $ follows plan $ p^* $
	}
}

\caption{\label{alg:insertion_heuristic} Insertion Heuristic}
\end{algorithm}

\subsubsection{Vehicle Group Assignment Method}
\label{sec:vga}
The VGA method relies on the performance improvement coming from conversion of a DARP problem to a variant of assignment problem.
In this work, we generalize the formulation by~\textcite{capMultiObjectiveAnalysisRidesharing2018} to be applicable in an online optimization setting. 
That is, we reformulated the algorithm to \rev{support optimization with requests already onboard some vehicles because the methodology by \textcite{capMultiObjectiveAnalysisRidesharing2018} assumes all vehicles to be empty before the request-vehicle matching.}

\begin{figure*}
\centering
\begin{subfigure}{0.9\columnwidth}
    \centering
    \includegraphics[width=.5\linewidth]{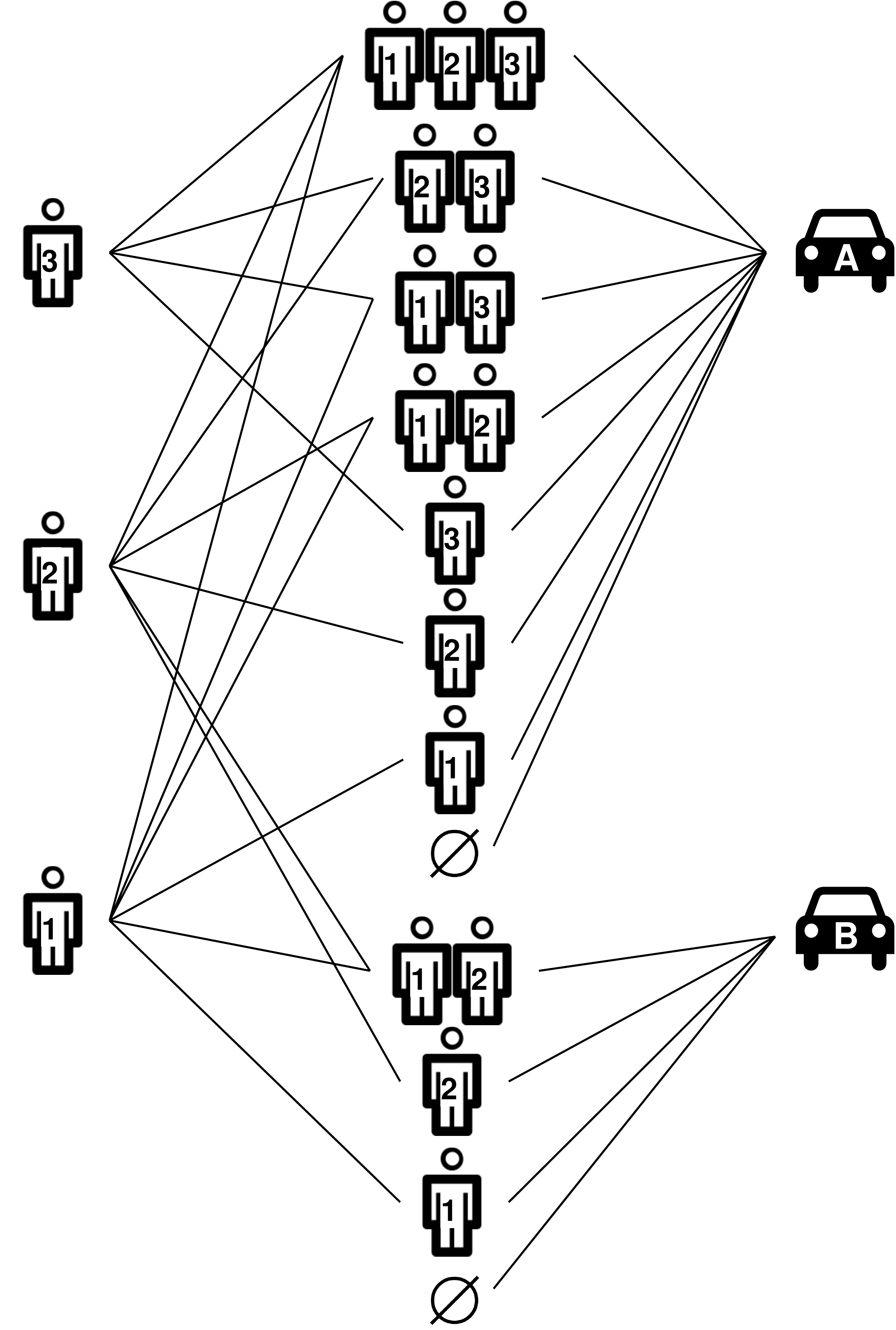}
    \caption{}\label{fig:vga_eample1}
\end{subfigure}
\begin{subfigure}{0.9\columnwidth}
    \centering
    \includegraphics[width=.5\linewidth]{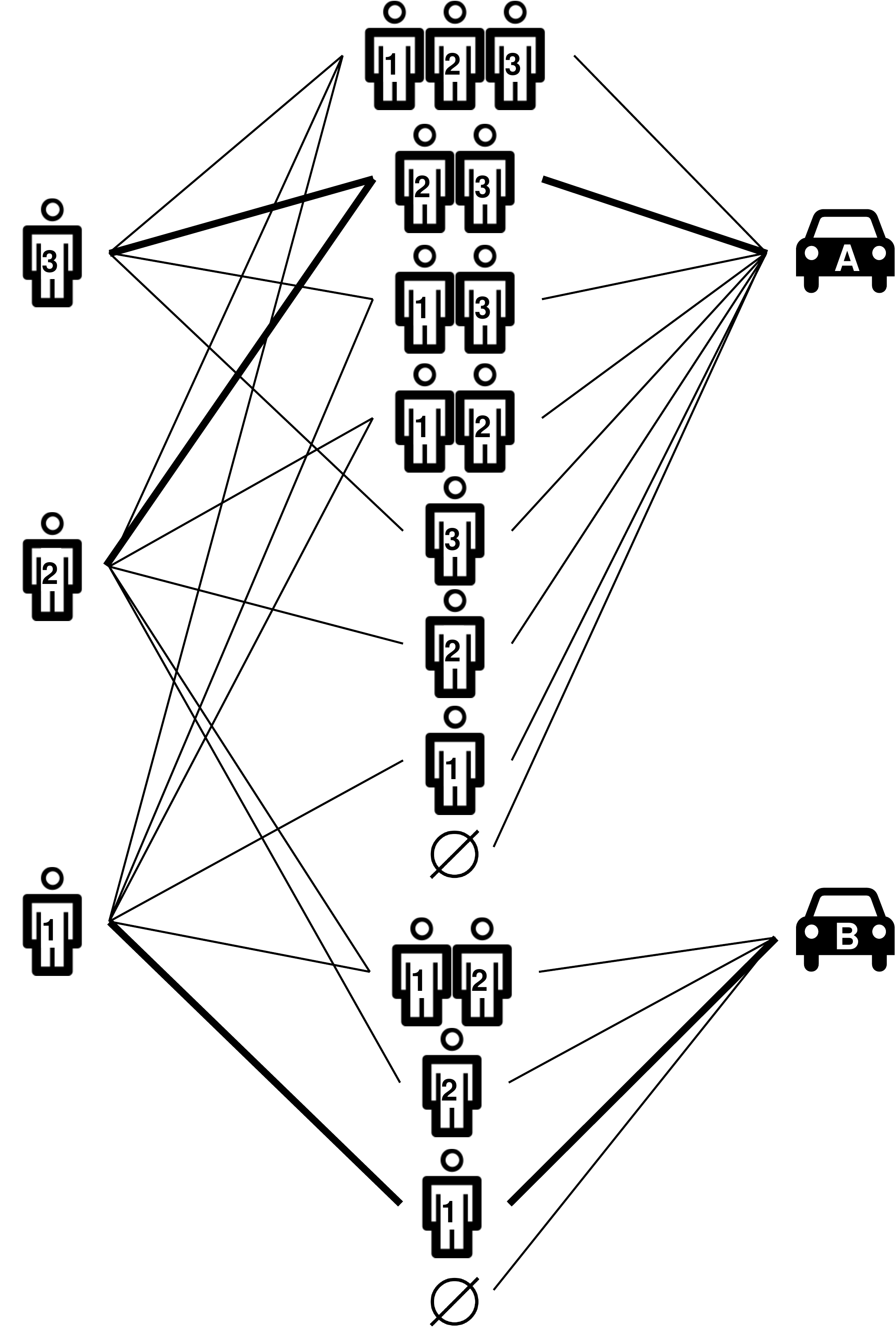}
    \caption{}\label{fig:vga_eample2}
\end{subfigure}
\caption{Example of the VGA method assigning three passengers to two vehicles.
In Figure~\ref{fig:vga_eample1}, we show all possible request groups for each vehicle. 
The lines between the request (left) and the group (middle) denote the membership in the group. 
The lines between the groups and vehicles denote feasible group assignments.
In Figure~\ref{fig:vga_eample2}, the final assignment between vehicles and groups is shown (bold lines).}
\label{fig:vga_example}
\end{figure*}

The VGA method can be divided into two phases: group generation (Algorithm~\ref{alg:dynamic-group-generation}) and vehicle-group assignment (Problem 1).
We can see the overall pseudocode in Algorithm~\ref{alg:vga}.
Let $ D_w $ be a set of waiting requests, i.e., the set of requests that have not been picked up yet.
Further, let \emph{group} be a set of requests such that for each group $ R $, $ R \subseteq D_w $.
In the first phase, for each vehicle, we compute all groups that can be serviced by the vehicle without violating the \rev{maximum delay $ q_{\mathrm{max}} $}
using the \emph{group generation} algorithm (Algorithm~\ref{alg:dynamic-group-generation}).
The second phase uses ILP (Problem~1) to map exactly one group to \rev{each} vehicle so that every request is serviced and the system plan is optimal.
The whole procedure is demonstrated by an example in Figure~\ref{fig:vga_example}.

\begin{algorithm}[t]
\Input{The current \rev{position and on-board passengers for each vehicle} in $ V $ and the set of waiting requests $ D_w $.}
\For{$ v \in V $}{
	$ \FG \leftarrow $ \texttt{generate\_groups($ v $, $ D_w $)}\;
}
$ \pi^* \leftarrow $ Solve Problem 1 using $ \FG[1] \dotsc \FG[m] $\;
All vehicles follow the optimal system plan $ \pi^* $\;
\caption{\label{alg:vga} VGA method}
\end{algorithm}

We say that a group $ R $ is \emph{feasible} for vehicle $ v $ if a plan exists for the vehicle that serves all requests from $ R $ \rev{without violating maximum delay and capacity constraints} and if all requests onboard vehicle $ v $ are members of the group.
We denote a set of all groups feasible for vehicle \( v \) by \( \FG \).
The key property of feasible groups, observed by~\textcite{alonso-moraOndemandHighcapacityRidesharing2017}, is that if a group \( R_1 \) is feasible, all subsets of the group \( R_2 \subset R_1 \) are also feasible.
This structural property is used to limit the number of groups we need to test for feasibility in the first part of the VGA method, the group generation algorithm.
To determine if a group is feasible, we define function \( \feas(R, v) \) that indicates whether the group \( R \) is feasible for vehicle \( v \).

The group generation algorithm (Algorithm~\ref{alg:dynamic-group-generation}) computes feasible groups for each vehicle independently.
First, the group generation algorithm computes all feasible requests for each vehicle $ v $, i.e., the feasible groups of size 1 marked as $ \FG[v][1] $ \rev{(lines 4-6)}.
Then, we find larger groups iteratively by combining the feasible groups from the previous iteration with all feasible requests \rev{(lines 7-14).
From the performance perspective, it is important to try each group only once (the purpose of \texttt{checked}) and also to check that all possible subsets of $ R $ of size $ |R| - 1 $ are present in $ \FG[v][|R| - 1] $ before checking group $ R $ for feasibility.
}
At the end of this step, we have a set of feasible groups for each vehicle, as it is illustrated in Figure~\ref{fig:vga_eample1}.

\begin{algorithm}[t]
\Input{A vehicle $ v $ and the set of waiting requests $ D_w $.}
\Output{Set of all feasible groups for the vehicle ($ \FG $).}
	\emph{Let $ R^{\mathrm{init}} $ be the set of requests onboard vehicle $ v $}\;
	$ k \leftarrow \max(|R^{\mathrm{init}}|,1) $\;
	$ \FG[v][k] \leftarrow \{R^{\mathrm{init}}\} $\;

	\For{$ R \in D_w $}{
		\If{$ \feas(\{R\}, v) $}{
			$ \FG[v][1] \leftarrow \{\{R\}\} \cup \FG[v][1] $\;
		}
	}
	\While{$ \FG[v][k] \neq \emptyset $}{
		$ \FG[v][k + 1] \leftarrow \emptyset $\;
		\tcc{not check groups repeatedly}
		$ \mathrm{checked} \leftarrow \emptyset $\;
		\ForAll{$ R \in \FG[v][k], \{ R \} \in \FG[v][1] $}{
			\If{$ (R \cup \{R\}) \notin \mathrm{checked} $ \And $ \forall R['] \subset (R \cup \{R\}), |R'| = k: R' \in \FG[v][k] $ \And $ \feas(R \cup \{R\}, v) $ }{
				 $ \FG[v][k + 1] \leftarrow (R \cup \{R\}) \cup \FG[v][k + 1] $\;
			}
			$ \mathrm{checked} \leftarrow \mathrm{checked} \cup (R \cup \{R\})$
		}
		$ k \leftarrow k + 1 $\;
	}
	\If{$ |R^{\mathrm{init}}| > 0 $}{
		$ \FG[v][1] \leftarrow \emptyset $\;
	}

	$ \FG \leftarrow \{\emptyset\} \cup \FG[v][1] \cup \FG[v][2] \cup \cdots \cup \FG[v][k] $\;

\caption{\label{alg:dynamic-group-generation}Function \texttt{generate\_groups} that generates groups for vehicle~$ v $. The boolean-valued function $ \feas(R, v) $ evaluates to true if vehicle~$ v $ can serve all requests from group~$ R $ \rev{without violating maximum delay and capacity constraints}.}
\end{algorithm}


The second part of the method finds the assignment of groups to vehicles that minimizes the total traveled distance resulting from vehicle plans such that for each vehicle, exactly one of the groups feasible for the vehicle is assigned, and all requests are served.
The assignment of groups to vehicles is formulated as an ILP.
There is a binary variable $ \xi_v^g $ for each possible vehicle-group assignment where $ \xi_v^g = 1 $ if a group $ g  \in \{ 1, \dotsc, |\FG| \} $ is assigned to vehicle $ v $ and $ \xi_v^g = 0 $ otherwise.
Using these variables, the problem is defined as:

\myprob{Vehicle-group Assignment}
\label{prob:vga}
\begin{equation*}
\min \sum_{v = 1}^{m} \sum_{g = 1}^{|\FG|}
\xi_v^g c(p_v^{g^*}),
\end{equation*}

subject to

\setcounter{equation}{0}
\begin{alignat}{2}
\sum_{g = 1}^{|\FG|} \xi_v^g &= 1 &\quad &\forall{v \in V} \label{vga-constr:1} \\
\sum_{v = 1}^{m} \sum_{\rev{R} = 1}^{|\FG|} \mathbf{1}_{\rev{R}}(r) \xi_v^{\rev{R}} &= 1 &\quad &\forall{R \in D_w} \label{vga-constr:2}
\end{alignat}

In the objective function, \( p_v^{g^*} \) denotes the optimal plan for vehicle $ v $ to serve group \( R_g \).
Constraint~\ref{vga-constr:1} states that only one group can be assigned to each vehicle.
Constraint~\ref{vga-constr:2} ensures that each request is served by exactly one vehicle plan. 
The indicator function $ \mathbf{1}_{R_g}(r) $ is equal to $ 1 $ if the request $ r $ is a member of the group $ R_g $ and $ 0 $ otherwise.

By solving the above described ILP, we obtain an optimal assignment of vehicles to feasible groups (see Figure~\ref{fig:vga_eample2} for example assignment).
This assignment can be directly translated into vehicle plans that replace vehicle plans from the previous iterations.

\subsubsection{\rev{Complexity, Optimality, and Implementation}}
\rev{
The worst-case complexity of computing an assignment of a set of waiting requests $ D_w $ to a set of vehicles $ V $ using IH is $ O(|D_w| \cdot |V| \cdot l) $ where $ l $ is the length of a plan and can be bounded as $l \leq K + |D_w| $.
For VGA, we analyze the complexity of each phase separately. 
For group generation, the computational complexity is dominated by the need to verify the feasibility of a group, represented by the call of function $ \feas(R, v) $.
Solving this function, in fact, equals solving a single-vehicle DARP (see Section~\ref{sec:vga_optimizatios} for more detail), which is an NP-hard~\cite{tothVehicleRoutingProblems2014} problem.
The vehicle-group assignment is then obtained by solving an ILP, which is also, in general, an NP-hard problem~\cite{schrijverTheoryLinearInteger1986}. 
Therefore, the VGA method can, in the worst case, require computational time that is exponential in the number of requests and vehicles. 
}
\trev{
However, DARP problem instances appearing in the context of large-scale MoD systems tend to have structural properties that are beneficial to the VGA algorithm. 
Specifically, since large-scale MoD systems are designed to provide quality of service comparable to using a private vehicle, the maximum waiting at pick-up is usually constrained to be less than a few minutes, and similarly limited is the maximum delay at destination.   
Such tight pick-up and drop-off time window constraints are used by the VGA algorithm to prune the feasible solution space. 
In practice, the maximum group size that requires a feasibility check tends to be relatively small, and also the total number of feasible groups tends to be within the grasp of existing ILP solvers. 
Under such conditions, the VGA algorithm is able to generate optimal results in practical computation time.}

\srev{
As for the optimality, IH is a heuristic approach, and as such, it cannot guarantee that the generated solution is optimal.
VGA method will generate an optimal solution if all feasible groups are generated exhaustively, and the ILP program is solved to optimality. 
For proof, see~\textcite{capMultiObjectiveAnalysisRidesharing2018}.
}

For our case study, we implemented the IH and the VGA method in Java. 
The ILP appearing in the VGA method is solved using Gurobi.\footnote{\url{http://www.gurobi.com/}}
The request-vehicle matching procedure is run every \SI{30}{\s} of the simulation for both IH and VGA.
For both methods, the maximum delay constraint is set to \SI{4}{\min}, and the vehicle capacity is set to five passengers.
The ILP solver in the VGA method computes the optimal solution with the maximum optimality gap of \SI{0.02}{\percent}.

\rev{
While the ability to compute optimal ridesharing assignments is essential to understand the limit of performance gains that can be achieved by ridesharing (and the gap between the optimal performance and the performance of the heuristic solutions), a practical deployment may impose constraints on the maximum run time of a ridesharing algorithm.
Therefore, we also tested a resource-constrained version of the VGA method with the ILP solver optimality gap set to \SI{0.5}{\percent} and group generation time-limited to \SI{60}{\ms} per vehicle.
We refer to this version as \emph{VGA limited}.
Also, for comparison, we reimplemented the method proposed by \textcite{alonso-moraOndemandHighcapacityRidesharing2017} in their Manhattan taxi ridesharing study.
As described in the supplemental material of~\cite{alonso-moraOndemandHighcapacityRidesharing2017}, this solution method employs specific heuristics and optimization cut-offs to achieve practical run time. 
We refer to this version as \emph{VGA PNAS}.
}

We compute passenger-vehicle assignments together with the simulation sequentially, and thus from a simulation perspective, the ridesharing computation is an instantaneous event.
In case of practical deployment, one could achieve sufficiently low wall-clock running time by computing on a computational cluster with many CPU cores because the VGA algorithm is easily parallelizable. 

\rev{
The existing variants of the VGA method~\cite{alonso-moraOndemandHighcapacityRidesharing2017,capMultiObjectiveAnalysisRidesharing2018} were not able to solve the ridesharing instances appearing in our case study to optimality within 24 hour limit on computational time or with a 60 GB memory limit. 
Therefore, we implemented performance optimization described in Appendix~\ref{sec:vga_optimizatios} that enabled us to find optimal solutions to our instances respecting these limits. 
}

\subsection{Simulation}
\label{sec:simulation}
In our experiments, we simulated the following five scenarios:
\begin{itemize}
\item \emph{Present state}: All the requests are served by private vehicles. 
The vehicles are parked at the request's start location, i.e.,  there is no delay. 
The number of used vehicles is equal to the number of requests, and the total distance traveled is equal to the sum of the shortest paths between the origins of all requests and their destinations.
\item \emph{MoD w/o ridesharing}: MoD system without ridesharing, the plans are computed using IH, and the vehicle capacity is set to one, i.e., the passengers are not allowed to share rides.
\item \emph{MoD w. IH Ridesharing}: MoD system with ridesharing computed by the IH.
\item \emph{MoD w. VGA Ridesharing (optimal)}: MoD system with ridesharing computed by the VGA method 
to optimality.
\item \emph{MoD w. VGA Ridesharing (runtime limited)}: MoD system with ridesharing  computed by the VGA method, with the group generation time-limited to \SI{60}{\ms} per vehicle and the ILP solver maximum optimality gap of \SI{0.5}{\percent}.
\item \rev{\emph{MoD w. VGA Ridesharing (PNAS)}: MoD system with ridesharing computed by the VGA method, with a set of timeouts/heuristics as described in~\textcite{alonso-moraOndemandHighcapacityRidesharing2017}}
\srev{, which we have reimplemented for this article.}
\end{itemize}

We simulate a morning peak time interval 7:00-8:00 and an off-peak time interval 11:00-12:00.
To avoid the ``cold start'' artifacts, the simulation begins 30 minutes before the analyzed time interval, at 6:30 and 10:30, respectively, but for subsequent analysis, we only use the data captured after the thirty-minute start period.
Including the \SI{30}{\minute} warm-up time, there are \num{122473} requests in the morning peak, and \num{42633} requests in the off-peak experiment.

\begin{figure*}[ht]
\centering
\begin{subfigure}{0.49\textwidth}
    \centering
    \includegraphics[width=1\linewidth]{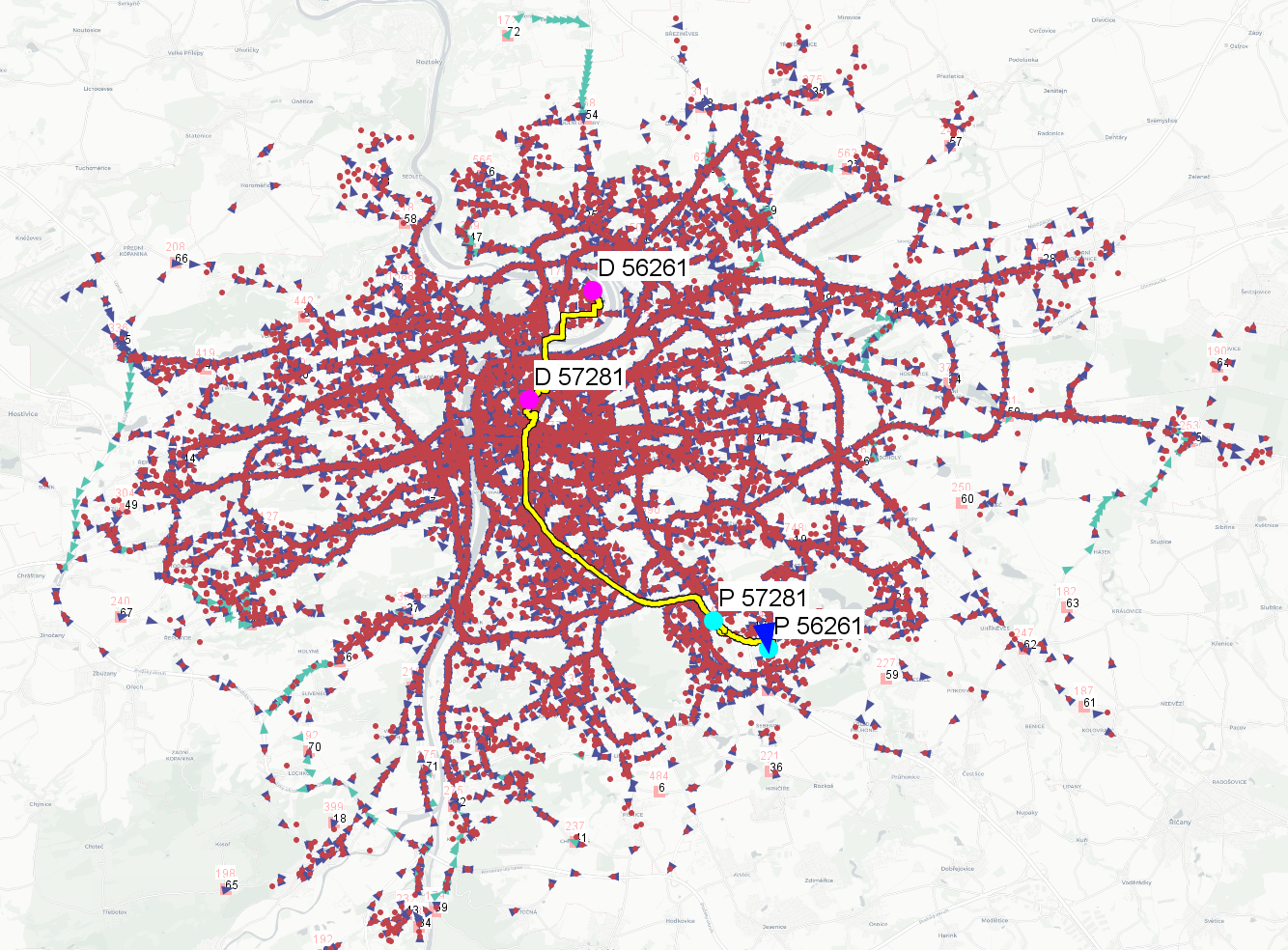}
    \caption{}\label{fig:ap-all}
\end{subfigure}
\begin{subfigure}{0.49\textwidth}
    \centering
    \includegraphics[width=1\linewidth]{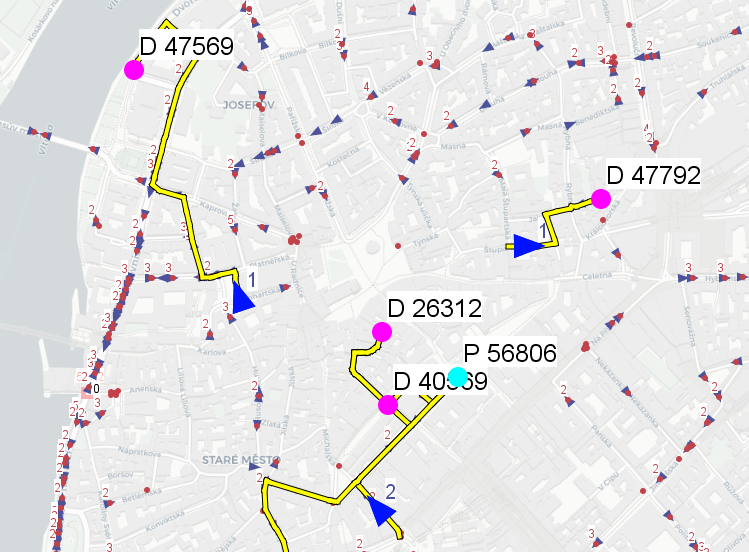}
    \caption{}\label{fig:ap-detail}
\end{subfigure}
\begin{tabular}{|c|}
\hline
\textbf{Video:} \url{https://sum.fel.cvut.cz/agentpolis/}\tabularnewline
\hline
\end{tabular}
\caption{\label{fig:ap}AgentPolis visualization of the simulated traffic in Prague during the traffic peak.
Figure~\ref{fig:ap-all} (left): the entire city of Prague in the simulation.
A more detailed (zoomed in) view can be seen in Figure~\ref{fig:ap-detail} (right).
Vehicles are represented as blue triangles, with a number indicating the onboard passenger count. Red circles represent passengers. 
Some vehicles are highlighted, and their current plan is drawn with a yellow line.
The pick-up and drop off locations of the remaining actions are marked with cyan and pink circles, respectively, with a number indicating passenger ID.
Note that in Figure~\ref{fig:ap-detail}, there are some passengers already driving in two of the vehicles, so the number of drop-off locations is greater than the number of pick-up locations.
The green triangles are vehicles that travel empty between stations (rebalancing).}
\end{figure*}

The scenarios were simulated in the multi-agent transportation simulation framework AgentPolis\footnote{\url{https://github.com/aicenter/agentpolis}}.
The simulation environment consists of a) road network composed of \emph{nodes} (crossroads) and \emph{edges} (road segments), b) on-demand vehicle stations, c) on-demand vehicle agents, and d) passenger agents.
In Figure~\ref{fig:ap}, we show a screenshot of the AgentPolis visualization captured during one of the simulation experiments.

We use an OpenStreetMap\footnote{\url{https://www.openstreetmap.org/}} road network consisting of \num{158674} edges and \num{63995} nodes. The speed limit for each road segment was also taken from OpenStreetMap data, and missing entries were generated according to following rules based on the local legislation: highway: \SI{130}{\km\per\hour}, living street: \SI{20}{\km\per\hour}, otherwise: \SI{50}{\km\per\hour}.

During initialization, we create vehicle stations, each filled with the pre-determined number of vehicles.
During the simulation, we are creating passenger agents for each request at its announcement time and origin point. 
Each passenger is then picked up by the assigned on-demand vehicle, driven to the desired location, dropped off, and finally released from the simulation.
The vehicle to serve the passenger is selected using the passenger-vehicle matching procedure (see Section~\ref{sec:matching}), either IH or VGA.
Note that each passenger can be either matched to one of the empty vehicles parked in a station or to a vehicle already serving some previously assigned requests.
Each vehicle executes its plan until it becomes empty (i.e., all assigned passengers have been dropped off), then it drives to park itself in the nearest station.
\rev{As explained in our assumptions, the passengers always select the MoD system as the mode for their trip.
In the simulation, any request that would wait for longer than 4 minutes or would be delivered to its destination with more than 4 minutes of delay is considered a rejected request. 
However, as mentioned before, we configured the system so that these quality of service bounds are always satisfied, and consequently, there are no rejections during the simulation experiment.}

\section{Results}
\label{sec:results}
In this section, we present the simulation results.
To run the experiments, we used a desktop system with Intel Core i7-8700K CPU~(\SI{3.7}{GHz}, 6/12 physical/virtual cores) and \SI{64}{GB} RAM.

\subsection{Operating Cost and Computational Time}
Tables~\ref{fig:comparison_table} and~\ref{fig:comparison_table-second_window} summarize the main results of the experiments. 
As explained in Section~\ref{sec:system_model}, we computed the size of the fleet to always guarantee full service availability. 
Since the service level is always \SI{100}{\percent}, we do not show this metric in result tables and plots.
The first row shows the value of our optimization criterion, i.e., the system operation cost measured in terms of total distance driven by the fleet vehicles. 
We can see that when using the VGA method instead of IH during the morning peak, we can save more than~\SI{110000}{\km} of vehicle distance driven, which represents more than~\SI{20}{\percent} reduction.
Compared to the ``no ridesharing" scenario and to the present state, the VGA method saves over \SI{573000}{\km} (\SI{57}{\percent}) and \SI{328000}{\km} (\SI{43}{\percent}), respectively.
Even in off-peak time, the VGA method can save about \SI{17}{\percent} of the total distance driven compared to the IH, and about \SI{48}{\percent} compared to the ``no ridesharing" scenario.

\begin{table}
\scriptsize
\centering{}%
\setlength{\tabcolsep}{0.3em}
{\renewcommand{\arraystretch}{1.2}%
\begin{tabular}{|+l|-r|-r|-r|-r|-r|-r|}
\hline
 &  & \multicolumn{5}{c|}{Mobility-on-Demand}
\tabularnewline
\cline{3-7}
 & \thead{Present} & \thead{No Ridesh.} & \thead{IH} & \thead{VGA} & \thead{VGA lim} & \thead{\rev{VGA PNAS}}
\tabularnewline
\hline
\hline
\rev{Optimal} & - & - & no & yes & no & no
\tabularnewline
\hline
\rowstyle{\bfseries}
Total veh. dist. (km) & \num{758001} & \num{1002766} & \num{539793} & \num{429172} & \rev{\num{451978}}& \num{475378}
\tabularnewline
\hline
Avg. delay (s) & - & \num{132} & \num{190} & \num{180} & \num{178} & \num{161}
\tabularnewline
\hline
\hline
Avg. density (veh/km) & \num{0.0077} & \num{0.0085} & \num{0.0053} & \num{0.0046} & \num{0.0048} & \num{0.0049}
\tabularnewline
\hline
Congested seg. & \num{8} & \num{25} & \num{1} & \num{1} & \num{1} & \num{0}
\tabularnewline
\hline
Heavily loaded seg.  & \num{163} & \num{291} & \num{31} & \num{10} & \num{17} & \num{20}
\tabularnewline
\hline
Used Vehicles  & \num{122473} & \num{33066} & \num{15685} & \num{13787} & \num{14449} & \num{14607}
\tabularnewline
\hline
Avg. comp. time (ms)  & - & \num{181} & \num{18} & \rev{\num{192903}} & \rev{\num{27714}} & \num{15598}
\tabularnewline
\hline
\end{tabular}}

\caption{\label{fig:comparison_table} Main results from the considered scenarios during the morning peak (7:00-8:00). Congested segments are segments on which traffic density is above critical density, and heavily loaded segments are segments with density above \SI{50}{\percent} of the critical density.}
\end{table}

\begin{table}
\scriptsize
\centering{}%
\setlength{\tabcolsep}{0.3em}
{\renewcommand{\arraystretch}{1.2}%
\begin{tabular}{|+l|-r|-r|-r|-r|-r|-r|}
\hline
 &  & \multicolumn{5}{c|}{Mobility-on-Demand}
\tabularnewline
\cline{3-7}
 & \thead{Present} & \thead{No Ridesh.} & \thead{IH} & \thead{VGA} & \thead{VGA lim} & \rev{\thead{VGA PNAS}}
\tabularnewline
\hline
\hline
Optimal & - & - & no & yes & no & no
\tabularnewline
\hline
\rowstyle{\bfseries}
Total veh. dist. (km) & \num{283483} & \num{344613} & \num{211285} & \num{175865} & \num{176957}& \num{186520}
\tabularnewline
\hline
Avg. delay (s) & - & \num{131} & \num{191} & \num{179} & \num{179} & \num{165}
\tabularnewline
\hline
\hline
Avg. density (veh/km) & \num{0.0045} & \num{0.0047} & \num{0.0034} & \num{0.0032} & \num{0.0032} & \num{0.0032}
\tabularnewline
\hline
Congested seg. & \num{0} & \num{1} & \num{1} & \num{0} & \num{0} & \num{0}
\tabularnewline
\hline
Heavily loaded seg.  & \num{5} & \num{10} & \num{3} & \num{1} & \num{2} & \num{2}
\tabularnewline
\hline
Used Vehicles  & \num{42633} & \num{7727} & \num{4646} & \num{4746} & \num{4802} & \num{5180}
\tabularnewline
\hline
Avg. comp. time (ms)  & - & \num{4} & \num{1} & \num{5438} & \num{4408} & \num{4113}
\tabularnewline
\hline
\end{tabular}}

\caption{\label{fig:comparison_table-second_window}Main results from the considered scenarios, off-peak (11:00-12:00). Congested segments are segments on which traffic density is above critical density, and heavily loaded segments are segments with density above \SI{50}{\percent} of the critical density.}
\end{table}

The VGA method is considerably slower than IH. The average computational time per one optimization batch in the peak scenario was about \SI{193}{\s}, compared to \SI{18}{\ms} for the IH.
Such a difference in the computational time may look extreme, but we have to consider the scale of the scenarios that were solved to optimality using the VGA method. 
The largest assignment problems (batches) contained more than \num{3000} waiting requests, \num{21000} active requests (including passengers already driving to their destination), and \num{11000} vehicles.

The runtime-limited experiment shows that we can speed up the VGA method significantly by merely limiting the computational time for the group generation and the solver.
In the VGA limited experiment, we reduce the computation time more than six-fold over the unconstrained version of the VGA method while still reducing the total traveled distance by more than \SI{16}{\percent} over the IH.
\rev{The VGA PNAS experiment reduces the computational time by another \SI{42}{\percent} at the cost of being closer to IH in the traveled distance (\SI{12}{\percent} improvement).
In the off-peak scenario, the VGA limited performs almost the same as the unconstrained version because the time limits are rarely reached.
Note, however, that there is more than a \SI{5}{\percent} increase of traveled distance in VGA PNAS, despite similar computational times suggesting that this method is not suitable for scenarios where sufficient computational resources are available.
}

\subsection{Trade-off Between Operating Cost and Passenger Discomfort}
\begin{figure}
\centering{}\includegraphics[width=1\columnwidth]{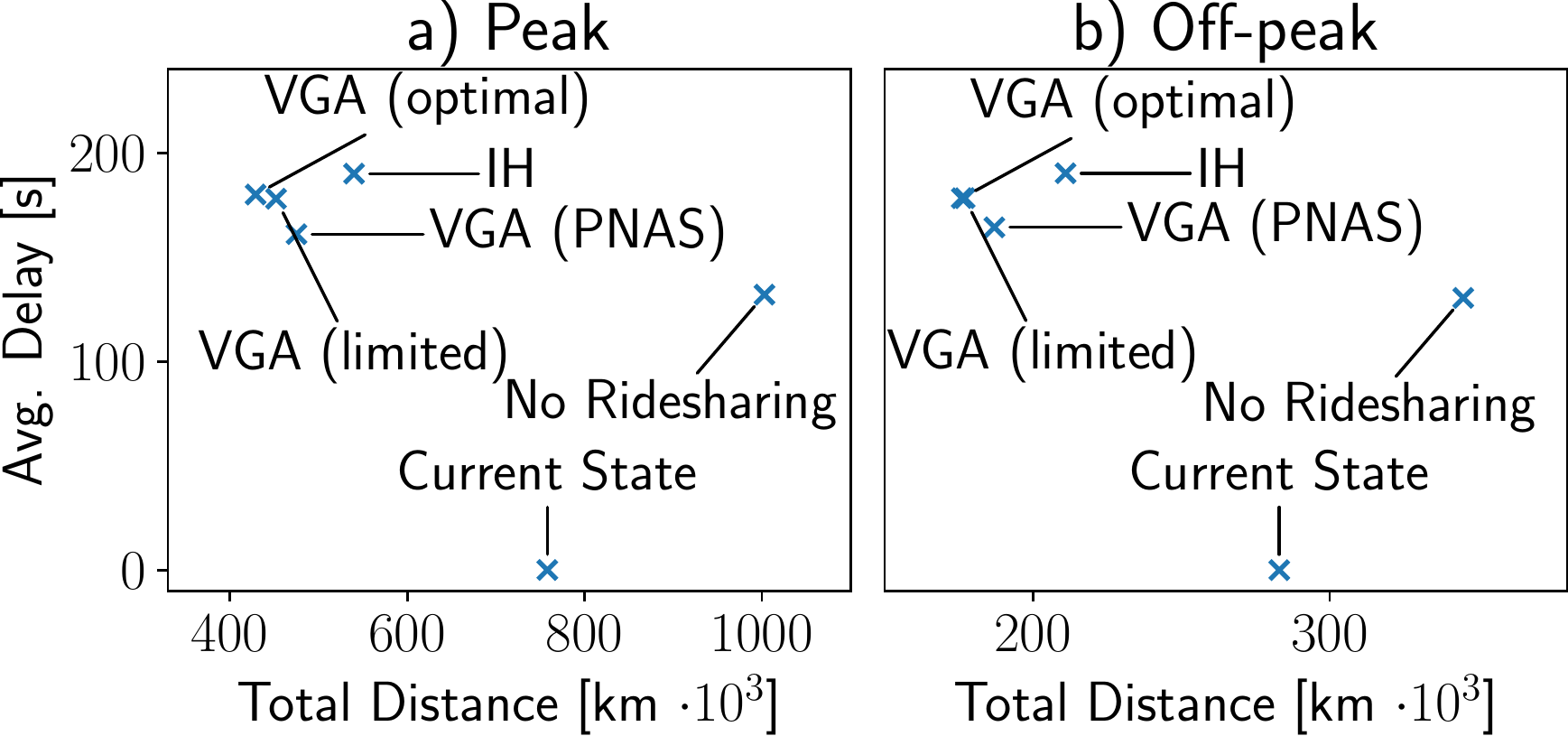}
\caption{\label{fig:distance_delay_tradeoff} The trade-off between total distance traveled by all vehicles and the average delay of one passenger trip.}
\end{figure}

Another metric that we tracked is the service quality, represented by the passenger delay relative to transportation by the private vehicle.
From Tables~\ref{fig:comparison_table}~and~\ref{fig:comparison_table-second_window}, we can see that the optimal VGA method saves about \SI{5}{\percent} time over the IH in both peak and off-peak experiments.
The trade-off between the operating cost (distance traveled) and the service quality (average delay) is depicted in Figure~\ref{fig:distance_delay_tradeoff}.

\begin{figure}
\centering{}\includegraphics[width=1\columnwidth]{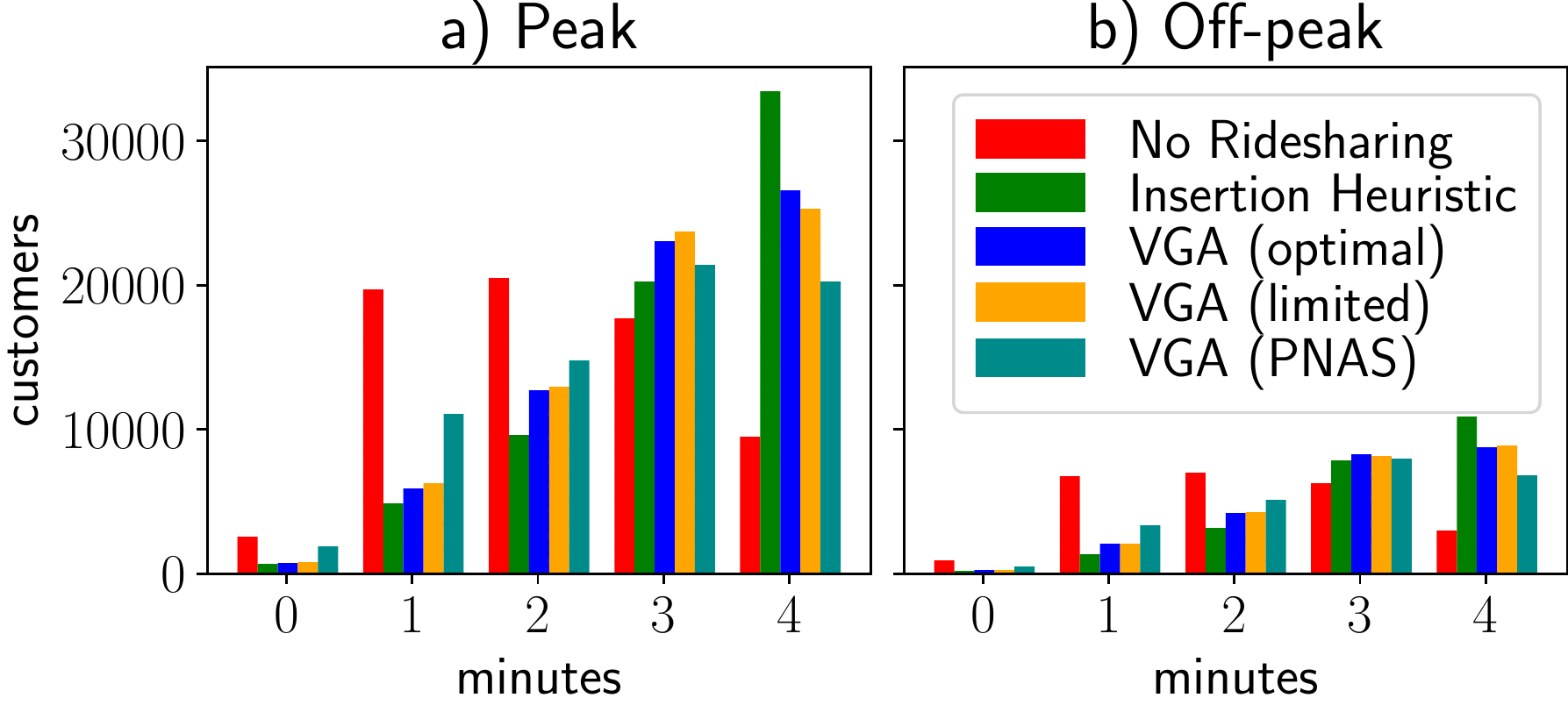}
\caption{\label{fig:delay_comparison} Histograms of delays. The present state scenario is omitted as the delay is always zero.}
\end{figure}

A more detailed overview of the passenger delays with a delay histogram for the four MoD scenarios in both time windows is in Figure~\ref{fig:delay_comparison}.
It is clear that for both peak and off-peak time, the VGA method reduces the passenger delay resulting from ridesharing compared to the IH.
Nevertheless, even in the case of the VGA method, there is a noticeably greater delay compared to the no ridesharing scenario, where the delay can occur only before the passenger is picked up or over the present state, where there is no delay because a car is assumed to be available at the origin of each passenger trip.

\subsection{Impact of MoD on Congestion}
\begin{figure*}
\centering{}\includegraphics[width=1\textwidth]{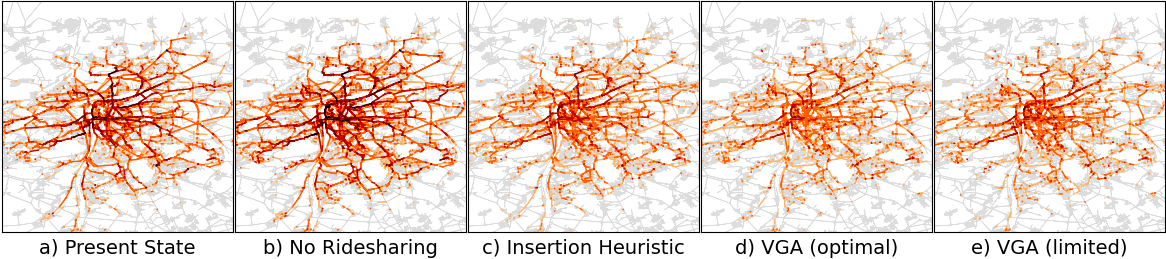}
\caption{\label{fig:traffic_density_map_comparison} Traffic density map of the four scenarios during the morning peak. Darker colors signalize higher traffic density. Black color means that the road segment is congested.
\rev{We omit the density map for the VGA PNAS experiment from this figure as it is very similar to the density map for the VGA limited experiment.}
}
\end{figure*}

\begin{figure*}
\centering{}\includegraphics[width=1\textwidth]{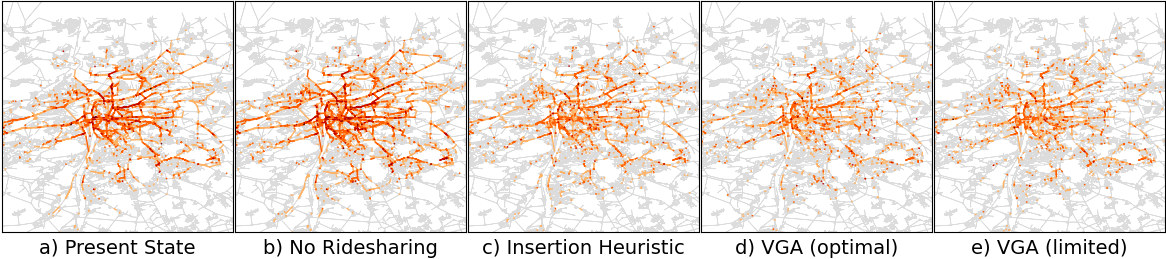}
\caption{\label{fig:traffic_density_map_comparison-second window} Traffic density map of the four scenarios during the off-peak time. Darker colors signalize higher traffic density. Black color means that the road segment is congested.
\rev{We omit the density map for the VGA PNAS experiment from this figure as it is very similar to the density map for the VGA limited experiment.}
}
\end{figure*}

In addition to the operational cost, we measured the impact of the MoD system on congestion.
We consider road segments with traffic density above the critical density of $ 0.08 \,\mathrm{vehicle}\, \mathrm{m}^{-1}$~\cite{tadakiCriticalDensityExperimental2015a}  as \emph{congested}. Segments with density above $ 0.04 \,\mathrm{vehicle}\, \mathrm{m}^{-1}$ are considered as \emph{heavily loaded}.
As you can see in Table~\ref{fig:comparison_table}, in the morning peak, using the optimal VGA method reduces the average traffic density by \SI{13}{\percent} over the ridesharing that uses IH, and by \SI{46}{\percent} and \SI{40}{\percent} over the MoD without ridesharing and the current state, respectively.
We can see the same trend when we look at the number of congested and heavily loaded segments.
In the off-peak experiment, the situation is similar, but the absolute numbers indeed show that there is no congestion in any of the scenarios.
Finally, Figures~\ref{fig:traffic_density_map_comparison} and~\ref{fig:traffic_density_map_comparison-second window} depict traffic densities on every road for all five scenarios.

\subsection{Fleet Size and Vehicle Occupancy}
Also, for each scenario, we recorded the number of vehicles that were used at least once during the simulation.
For the present state scenario, we consider a dedicated vehicle for each request. Therefore, the number of used vehicles is equal to the number of requests. 
The results confirm that the VGA method indeed makes the MoD system more efficient.
During peak-hour, the optimal VGA used \rev{1898} (\SI{12}{\percent}) fewer vehicles than the IH.  
Compared to the MoD system without ridesharing, the MoD system with optimal ridesharing used about one-third of the vehicle fleet, and compared to the present state system, the reduction is almost thirteen-fold.

\begin{figure*}
\centering
\begin{subfigure}{0.245\linewidth}
    \centering
    \includegraphics[width=1\linewidth]{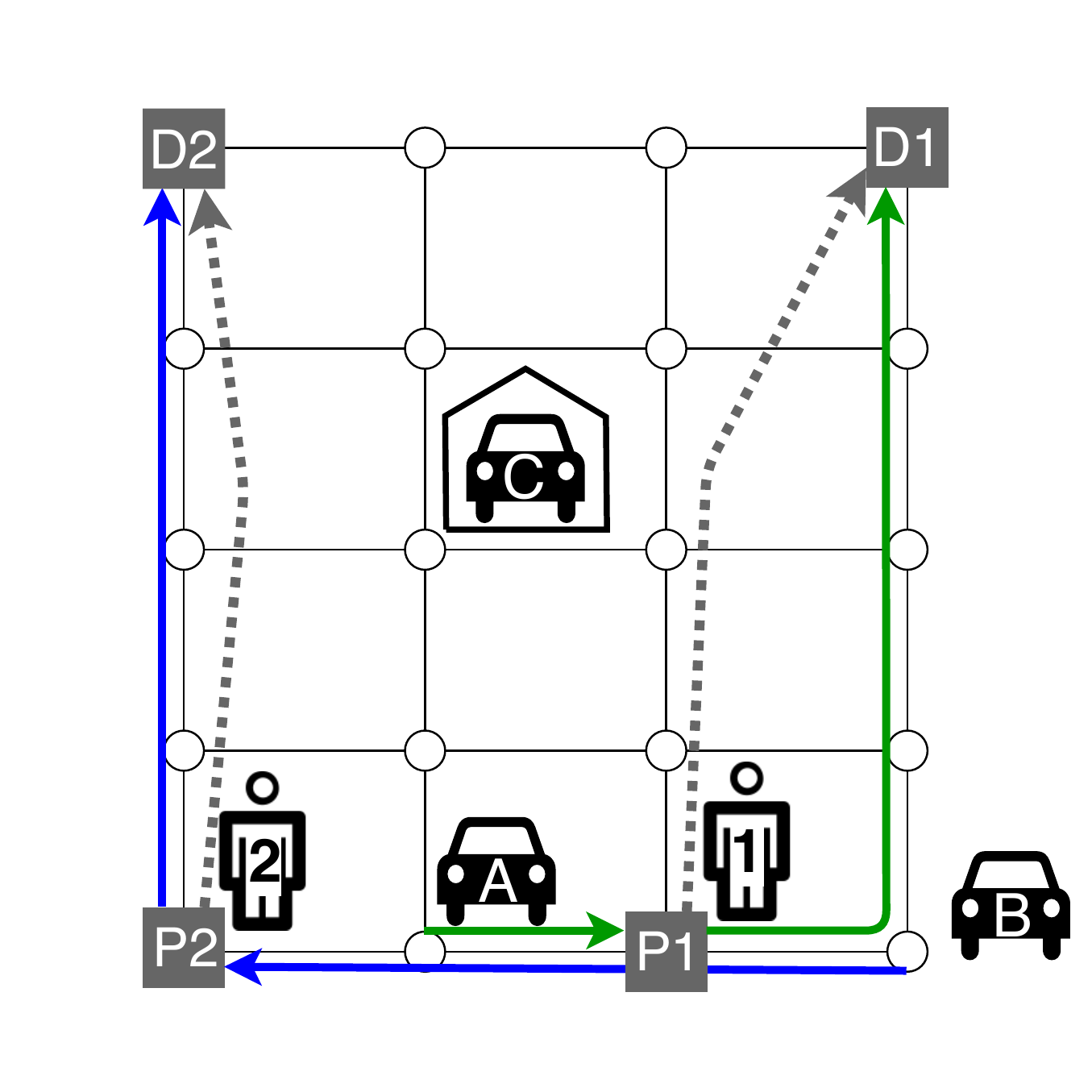}
    \caption{IH Iteration 1}\label{fig:car_example-IH-a}
\end{subfigure}
\begin{subfigure}{0.245\linewidth}
    \centering
    \includegraphics[width=1\linewidth]{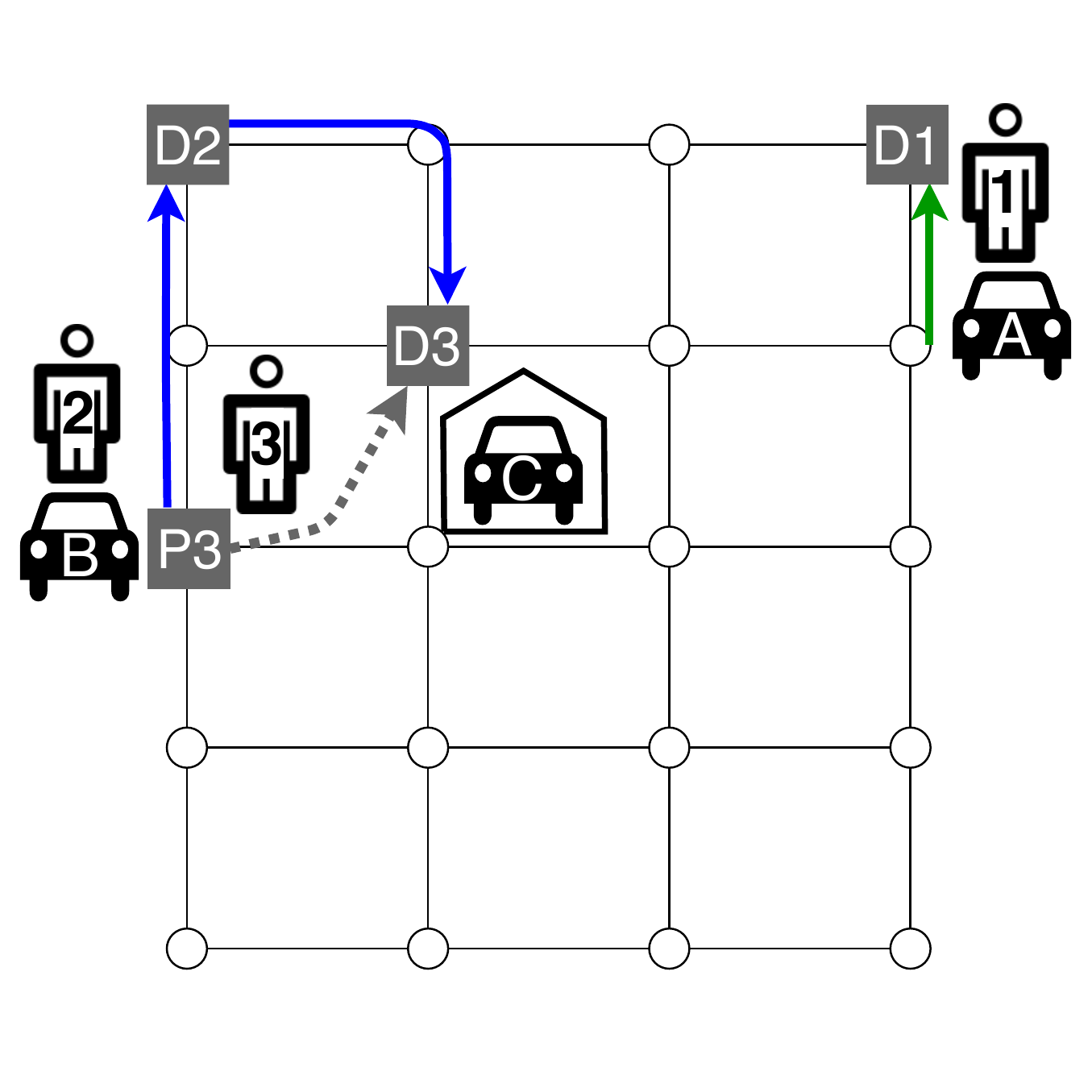}
    \caption{IH Iteration 2}\label{fig:car_example-IH-b}
\end{subfigure}
\begin{subfigure}{0.245\linewidth}
    \centering
    \includegraphics[width=1\linewidth]{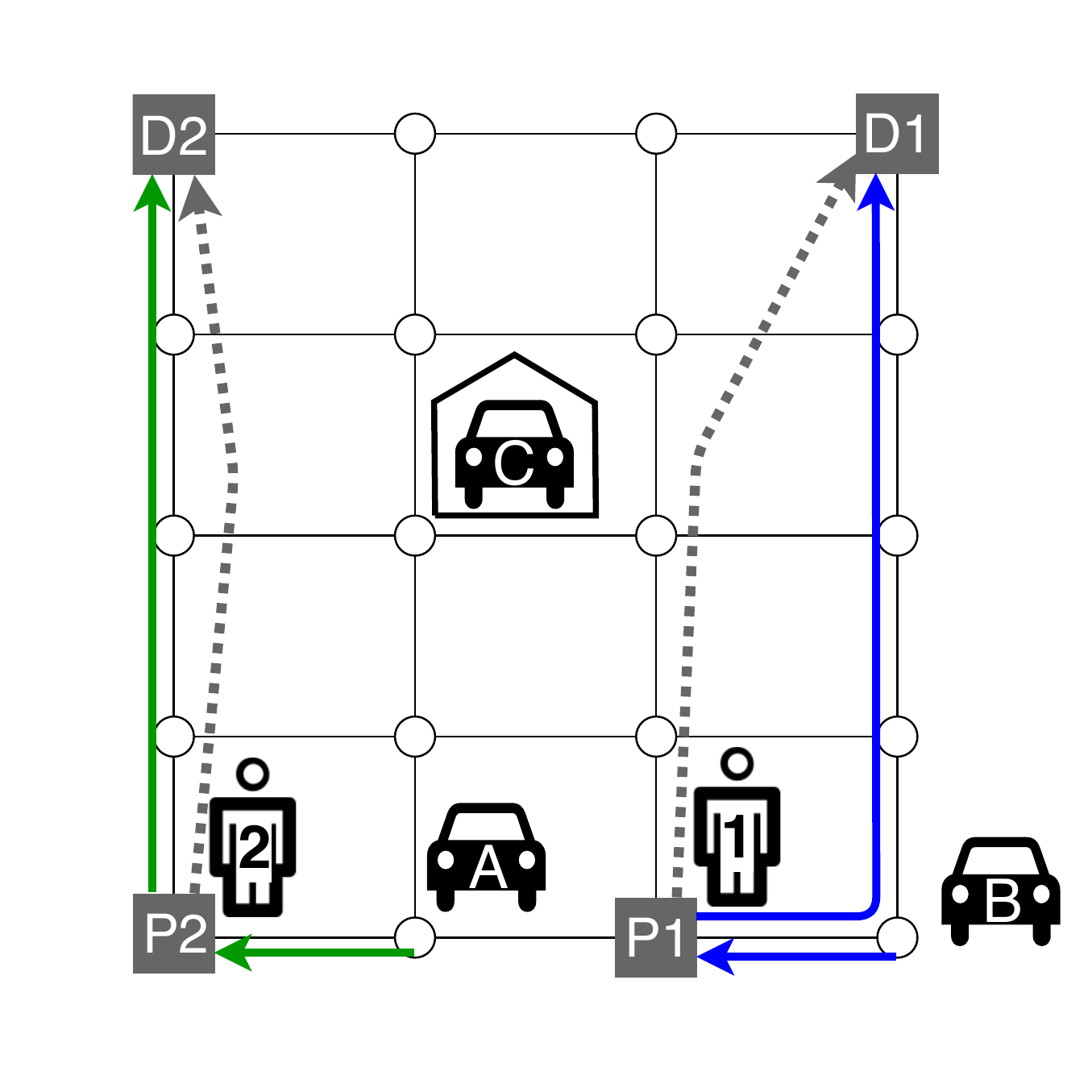}
    \caption{VGA Iteration 1}\label{fig:car_example-VGA-a}
\end{subfigure}
\begin{subfigure}{0.245\linewidth}
    \centering
    \includegraphics[width=1\linewidth]{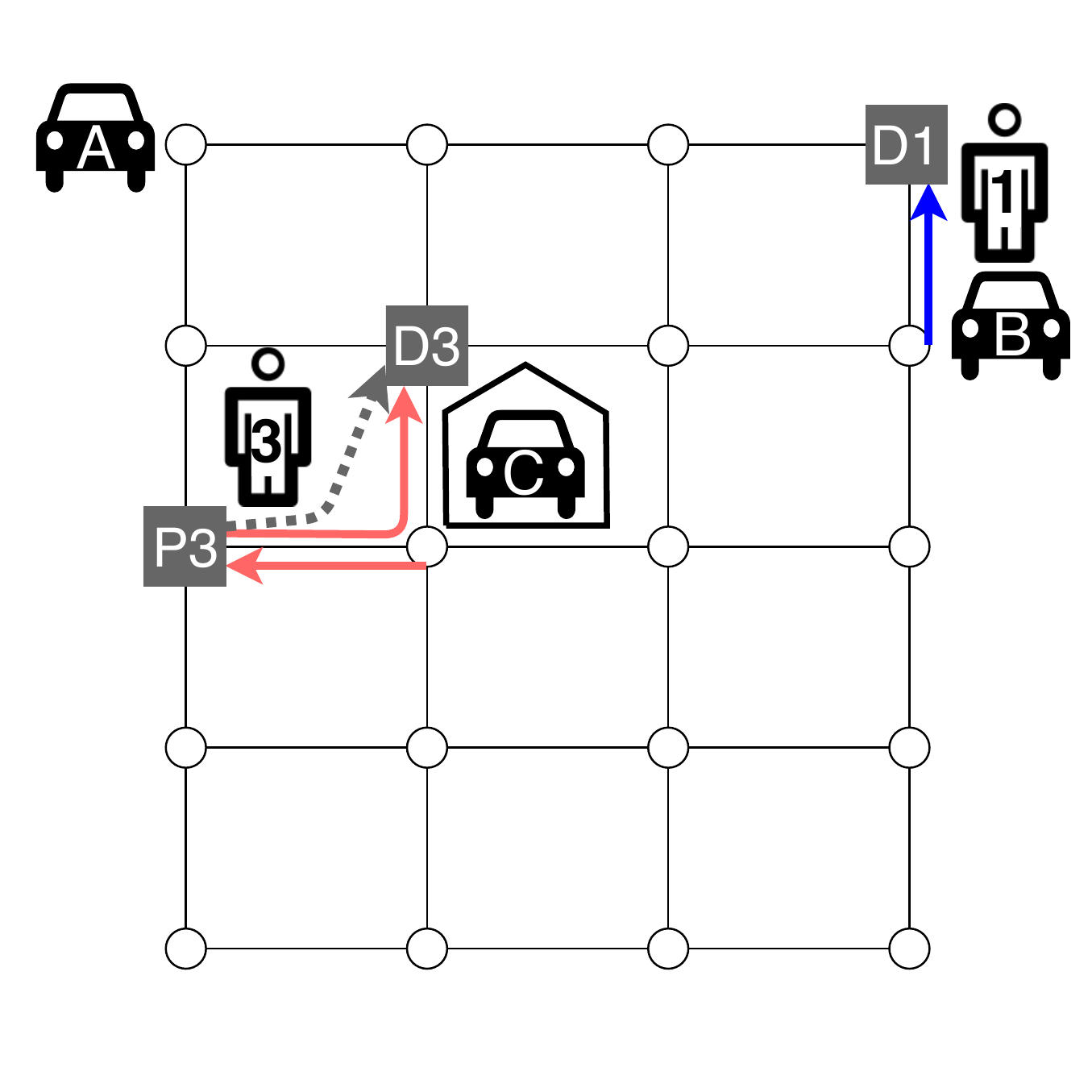}
    \caption{VGA Iteration 2}\label{fig:car_example-VGA-b}
\end{subfigure}
\caption{Example of the capital cost paradox.
In Figures~\ref{fig:car_example-IH-a}~and~\ref{fig:car_example-IH-b}, we can see two iterations of the IH.
In Figure~\ref{fig:car_example-IH-a}, there are three vehicles: vehicles $ A $ and $ B $, and vehicle $ C $ that resides in the station, representing a potentially unlimited pool of vehicles.
Also, there are two passengers ($ 1 $ and $ 2 $), that request the travel from their current locations  $ P1 $ and $ P2 $ to their destinations $ D1 $ and $ D2 $ (denoted by dashed arrows). 
Solid arrows denote the plans for both vehicles computed by the first iteration of the IH.
In Figure~\ref{fig:car_example-IH-b}, there is the same scenario in the next iteration. 
Both cars moved by five steps in the grid, and also, a new request appeared. 
We can see the new plans generated by the second iteration of IH too.
The second set of Figures~((\ref{fig:car_example-VGA-a})~and~(\ref{fig:car_example-VGA-b})) shows the exact same two iterations solved by the VGA method.
Note that although we saved one segment of traveled distance (vehicles traveled 14 segments in the grid combined compared to 15 segments in case of the IH), we used one extra vehicle (vehicle $ C $) that was not needed in the IH scenario, thus effectively increased the required fleet. 
}
\label{fig:paradox}
\end{figure*}

In the off-peak time, however, we registered that the optimal VGA method uses about \SI{2}{\percent} more vehicles than IH.
By analyzing the simulation output, we found an explanation for this perhaps surprising result.
First, counterintuitively, it is possible that a suboptimal vehicle assignment that generates plans with longer total distance can lead to fewer vehicles being used, as it is illustrated in Figure~\ref{fig:paradox}.
Second, by analyzing the vehicle trips in both IH and VGA scenario, we found that such situations occur frequently due to unbalanced demand.
In other words, the optimal method uses more vehicles not despite, but because its plans are more operating cost-efficient: the vehicles simply serve requests too quickly, which increase the chance of ending up in the areas with lower demand, where they need to wait a long time before another request appears nearby.
This reminds us that to fully understand MoD systems, we need to study not only operation cost vs. service quality trade-offs, but also operation-cost vs. capital cost trade-offs associated with different design and control strategies. 

\begin{figure}
\hspace*{-1cm}
\centering{}\includegraphics[width=1.1\columnwidth]{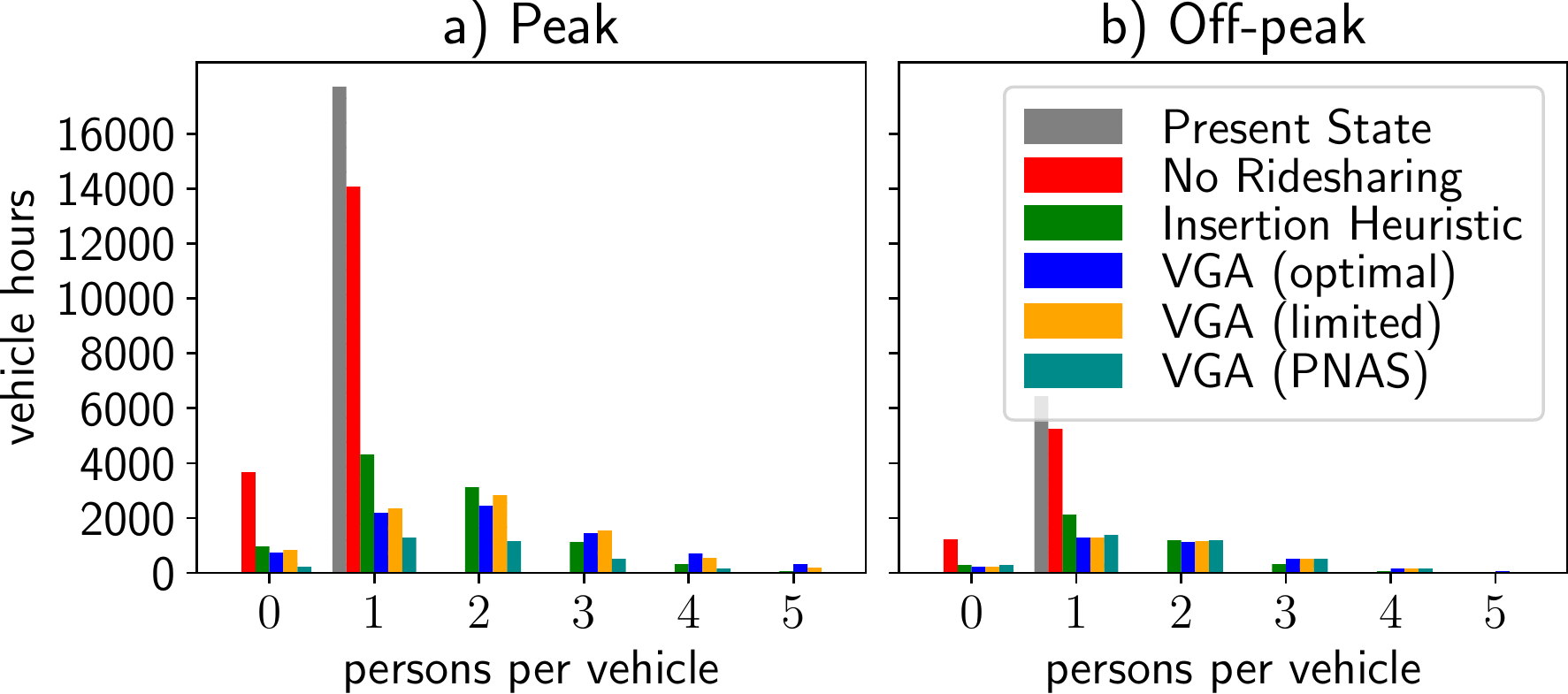}
\caption{\label{fig:occupancy_comparison} Occupancy histogram of all five scenarios.}
\end{figure}

Next, we measured vehicle occupancy: Figure~\ref{fig:occupancy_comparison} shows the occupancy histogram for the four compared scenarios.
We can see that vehicle occupancy is the highest when using the optimal method in both peak and off-peak scenarios.

\subsection{\rev{Sensitivity Analysis of the VGA Method}}
\begin{figure}
\centering{}\includegraphics[width=0.9\columnwidth]{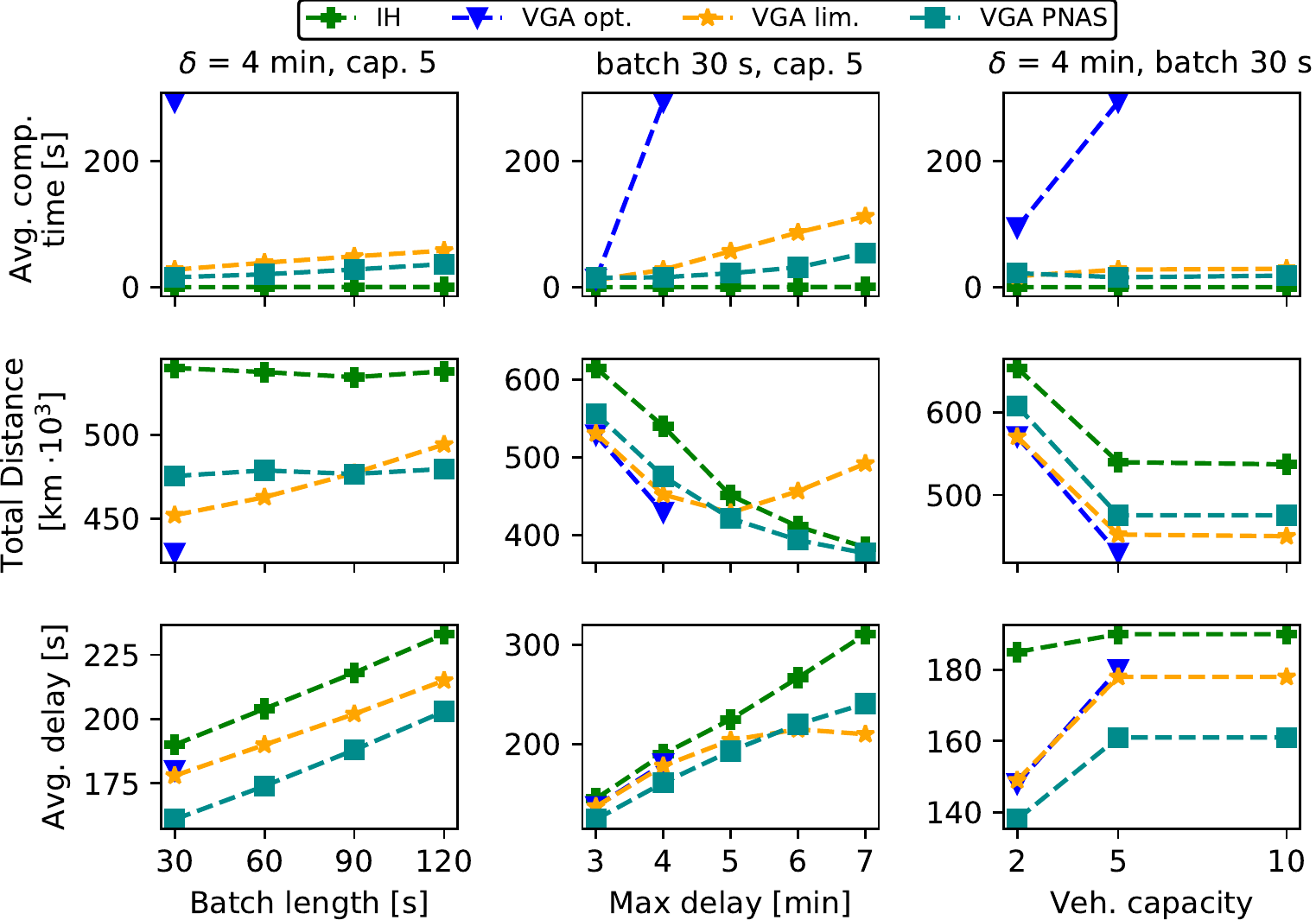}
\caption{\label{fig:sensitivity_analysis-peak} \srev{Sensitivity analysis: peak. 
Each column represents one experiment set, and inside each column, each value on the x-axis represents one experiment. 
Each row displays a single measured quantity. 
The optimal VGA method is only computed for batch length \SI{30}{\s}, maximum delay 3 and 4 minutes, and capacity 2 and 5 persons per vehicle. 
For other parameter values, the optimal VGA method did not terminate within 24 hours.}}
\end{figure}

\begin{figure}
\centering{}\includegraphics[width=0.9\columnwidth]{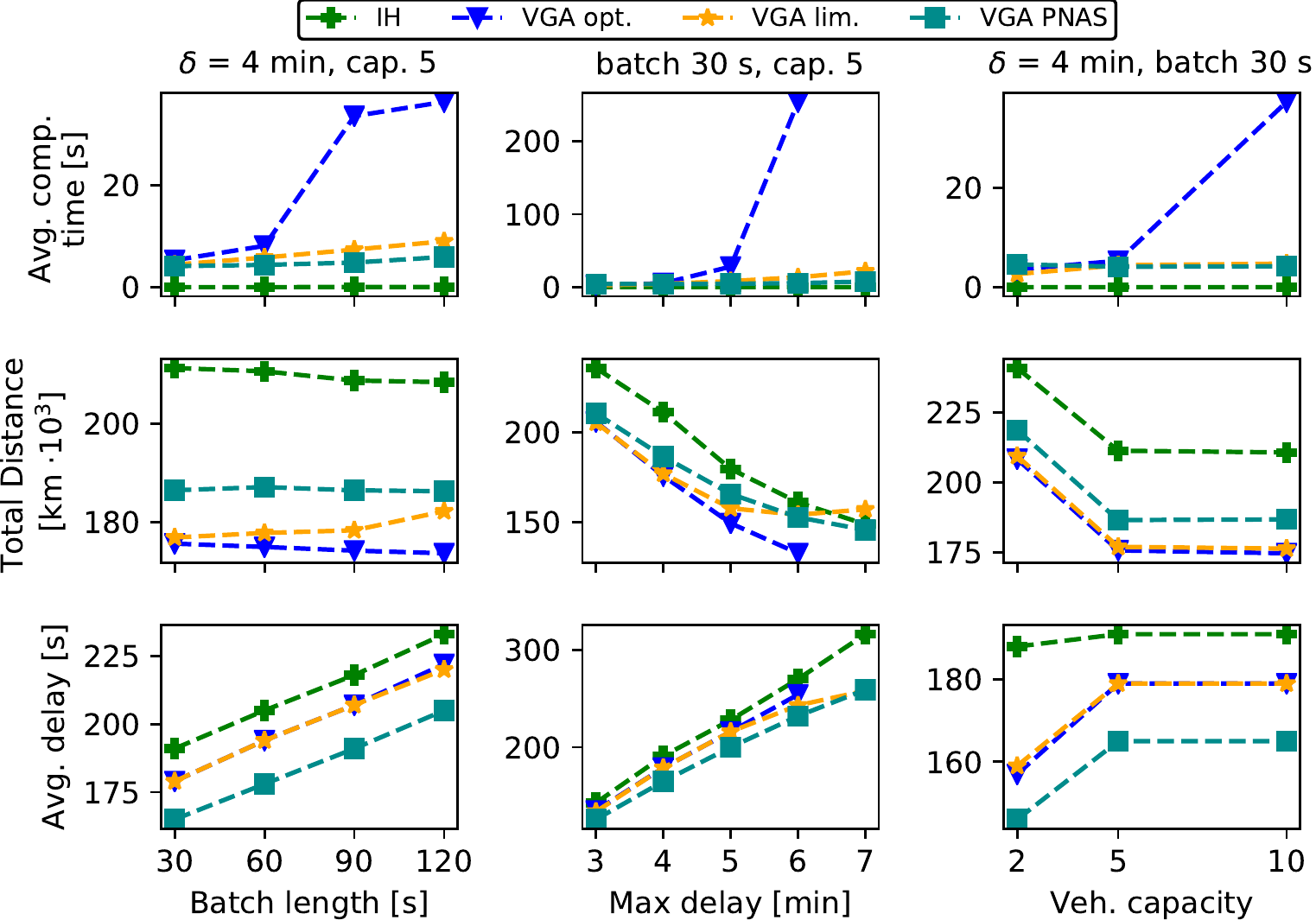}
\caption{\label{fig:sensitivity_analysis} \srev{Sensitivity analysis: off-peak. 
Each column represents one experiment set, and inside each column, each value on the x-axis represents one experiment. 
Each row displays a single measured quantity. 
The optimal VGA method is only computed for the maximum delay of up to 6 minutes.
For the maximum delay of 7 minutes, the optimal ridesharing assignment cannot be computed within 24 hours limit.}}
\end{figure}

\rev{
We analyze the sensitivity of \srev{IH and} the three variants of the VGA method (optimal, limited, PNAS) to variation in batch length, maximum delay\srev{, and capacity} with respect to total traveled distance, computation time, and average passenger delay.
Note that the time between a request announcement and the end of the batch, when the passenger-vehicle assignment is recomputed, counts towards the delay of the request.
Therefore, for scenarios with longer batches, we also extended the maximum delay in order to keep the average effective maximum delay of \SI{4}{minutes}.
Also, note that we use the same stations and fleet for all experiments, and consequently, some requests were rejected in configurations with shorter maximum delay or longer batch length. 
However, the service level remains above \SI{99}{\percent} in all configurations, so the impact of rejected requests on the results is negligible.

We can see the results in Figure~\ref{fig:sensitivity_analysis-peak}~(peak scenario) and Figure~\ref{fig:sensitivity_analysis} (off-peak scenario).
In the peak scenario, we were able to compute the optimal solution only for batch length of up to 30 seconds and for up maximum delay of up to \SI{4}{minutes}. 
For larger values, the algorithm failed to terminate within \SI{24}{hours}.
As expected, the runtime of the optimal method grows exponentially with maximum delay.
This is best seen in the case of the off-peak experiment\srev{, where the algorithm was able to find an optimal solution within 5 minutes of runtime on average for the 6-minute maximum delay but failed to compute optimal solutions within 24 hours runtime for the 7-minute maximum delay.}
Clearly, the optimal method  would not scale to scenarios with larger limits on maximum delay, and one of the resource-limited variants would have to be employed.

\srev{We can observe that the system's efficiency (measured in terms of total distance driven) monotonically increases with maximum allowed delay for all considered methods with the exception of the VGA limited method, which achieves low runtimes by prematurely terminating computation in several stages of the algorithm. 
As we can see, the values of cut-off parameters that work well for 3-5 minute maximum delay lead to inferior performance for larger maximum delays.
We can also observe that with the increasing maximum delay, the gap between the optimal method and both resource-constrained VGA methods increases.
Remarkably, our experiments show that for high values of max delay, the IH achieves almost identical performance as the VGA PNAS algorithm while using only a fraction of computational resources.  
}


The batch length negatively impacts the average travel delay experienced by passengers because they need to wait for the end of the batch for their request to be assigned to a vehicle. 
The motivation for using longer batch lengths is to gather more requests and to find a more efficient passenger-vehicle assignment.
However, it appears that even in the off-peak experiment, these efficiency gains get only realized using the optimal solution method. 
For suboptimal solution methods, the efficiency gains are either negligible or straight-out negative.
\srev{The higher vehicle capacity also increases the potential for ridesharing, which, in turn, could improve the operational efficiency of the system.
In the peak experiment, the increased complexity prevented the computation of the optimal solution for vehicle capacity set to 10 passengers.
However, our results show that the resource-limited variants were not able to improve the solution significantly. 
This can be caused by reaching the time limits before the algorithm can generate high-occupancy plans, or it can indicate that the capacity of 5 passengers per vehicle is sufficient.
For reference, in the off-peak experiment, we were able to compute optimal solutions, and our results show that the vehicle capacity of 5 is sufficient for the off-peak demand intensity. However, it is still possible that high-capacity vehicles would be better utilized when planning for peak-intensity demand.
}

Concerning the suboptimal versions of the VGA method, our VGA limited method computes higher quality solutions for shorter batch lengths and shorter maximum travel delays, while the PNAS version of the VGA method achieves better performance for the longer batch lengths and longer maximum delays. 
This is probably because the PNAS version uses IH to compute plans for groups larger than four (see the supplemental material of~\textcite{alonso-moraOndemandHighcapacityRidesharing2017}).
This will negatively affect solution quality for easier instances. 
However, for harder instances, this approach may be beneficial compared to our strategy because it allows forming larger groups with potentially suboptimal plans.
This observation suggests that a resource-constrained VGA method should use a heuristic to compute larger groups, but the threshold for using this heuristic should be determined by the remaining computational time.
}

\subsection{Computational Time Analysis of the VGA Method}
\label{sec:performance}
\begin{figure}[ht]
\centering{}\includegraphics[width=1.0\columnwidth]{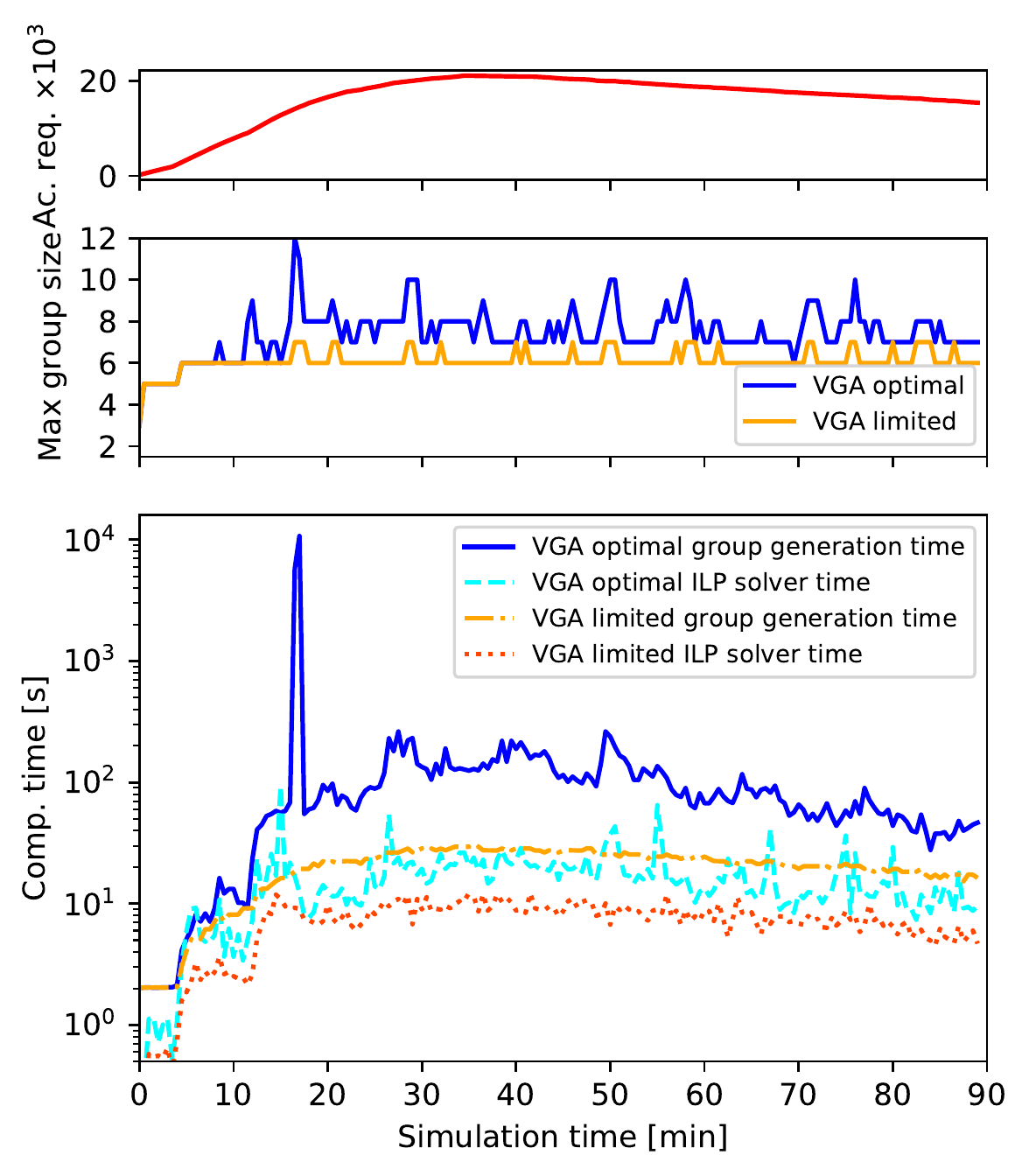}
\caption{\label{fig:vga-simulation-time} Computational efficiency analysis of the VGA scenarios during the peak time.
In the top figure, we show the evolution of the number of active requests over time. 
We can see that after the warm-up time, the number of active requests in the system is stable, only slowly decreasing.
The middle figure displays the maximum group size that was computed in each batch.
The bottom figure demonstrates how the computational time, consisting of the group generation time and the ILP solver time, change during the simulation.
}
\end{figure}

Finally, we inspect the computational requirements of the VGA method.
In Figure~\ref{fig:vga-simulation-time}, we show the evolution of the number of active requests (top), maximum computed group size (middle), and computational time for group generation and group-vehicle assignment process (bottom) during the peak scenarios, including the warm-up period.
Looking at the maximum group size, we see that the \SI{60}{\ms} limit for the group generation results in groups of the maximum size of 5-7 in most batches, while in the optimal scenario, the maximum group size has high variance and goes up to \num{11}. 

When we compare the maximum group size with the computation times, we can obtain other valuable insights: a) the group generation time is strongly dependent on the maximum group size, and thus it has low variance in the limited scenario and high variation in the optimal scenario, b) the solver time does not depend on maximum group generation time much, and it is highly variable in both limited and optimal variant, and c) the group generation time dominates \rev{in both scenarios}.
\rev{These findings suggest that further performance optimization of the group generation process may lead to a more favorable trade-off between the solution cost and computational time.}

\section{Conclusion}
\label{sec:conclusion}

Urban MoD systems represent a promising alternative to private car transport that can reduce the number of vehicles by employing massive vehicle sharing.
To further improve the efficiency of an MoD system, the system operator can implement large-scale ridesharing, where multiple passengers are transported in one vehicle simultaneously.
Ridesharing can increase vehicle occupancy and reduce the total distance driven in the system,
but finding the optimal assignment of passengers to vehicles is a hard combinatorial problem.
Traditional exact algorithms for vehicle routing are only applicable to the instances that are orders of magnitude smaller than instances occurring in the metropolitan-scale MoD systems.
Therefore simpler heuristic methods for ridesharing are often employed.
Recently, the Vehicle-Group Assignment (VGA) has been shown to be capable of solving ridesharing problems with up to \num{500} vehicles and requests optimally. 

In this work, we implemented algorithmic improvements that allowed us to successfully apply the VGA method to a metropolitan-scale MoD system.
In contrast to previous studies that sacrifice either scale or optimality, we can regularly compute optimal assignments of more than \num{21000} active requests to over \num{10000} vehicles. 
Also, we study the trade-off between the MoD system efficiency and computational performance for several other passenger-vehicle assignment methods.
Specifically, we compared \rev{six} different scenarios: 1) the "status quo" system with private vehicles, 2) MoD system without ridesharing, 3) MoD system with ridesharing using IH, 4) MoD system with ridesharing using optimal assignments computed by the VGA method, \rev{5, 6) MoD systems with ridesharing that use two resource-limited versions of the VGA method.}
For all five scenarios, we measured operation cost (total vehicle distance driven), service quality (average delay), fleet size, and congestion levels. 
Also, we measured the computational time for ridesharing methods, and in the case of the VGA method, we performed an analysis of the contribution of different sub-problems to the overall computational time.

The results confirmed that ridesharing dramatically increases the efficiency of an MoD system: by employing the VGA method, we reduced the total distance driven in the system by more than \SI{57}{\percent} compared to the present state.
Moreover, we demonstrated that the optimal ridesharing assignments are significantly more efficient than assignments computed by the heuristic approach.
Our results show that by using the optimal method instead of the IH, we can reduce the total distance traveled by more than \SI{20}{\percent} while simultaneously reducing the average passenger delay by \SI{5}{\percent}.
Finally, our resource-constrained VGA method provides more than \SI{16}{\percent} travel distance saving over IH while reducing the computational time by almost \SI{90}{\percent}.
\rev{Besides the expected conclusion that ridesharing yields significant savings, these results identify and quantify the optimality gap between a previously proposed resource-constrained version of the VGA method and the optimal solution, and also between a resource constrained VGA method and an IH-based ridesharing method. 
\frev{
Our sensitivity analysis provided insights into the limits of the VGA method. The method is capable of finding optimal assignments given that vehicle capacity is no more than 5-10 passengers and the maximum allowed delay is no more than 4-7 minutes, depending on the intensity and structure of the demand. However, for scenarios using higher-capacity vehicles or with more permissive delay constraints, the VGA algorithm can no longer certify the optimality of the computed ridesharing assignments. 

}
We believe that all these findings can help future researchers and practitioners to understand the trade-offs between different MoD system operating policies.
Moreover, the optimality gaps provide insights into maximum efficiency gain that one can hope to achieve by developing new heuristic solutions.}

\trev{
In future work, we plan to include more advanced metaheuristics in the comparison.
Also, we plan to consider a more general model of a mobility-on-demand system, where travellers could transfer between different vehicles, some of them potentially being fixed-route high-capacity vehicles such as buses or trains.
Finally, we plan to investigate the process of MoD system design, including fleet-sizing, fleet composition, and MoD operation from a multi-objective perspective, studying trade-offs between capital cost, operation cost, and service quality.
}

\section*{Funding}
This work was supported by the Czech Science Foundation under Grant No. 18-23623S; AMS Institute, and OP VVV MEYS funded project under Grant No. CZ.02.1.01/0.0/0.0/16\_019/0000765 ``Research Center for Informatics". Access to computing and storage facilities owned by parties and projects contributing to the National Grid Infrastructure MetaCentrum, provided under the program ``Projects of Large Infrastructure for Research, Development, and Innovations" (LM2010005), is greatly appreciated.


\printbibliography

@article{pavone_robotic_2012,
  title = {Robotic Load Balancing for Mobility-on-Demand Systems},
  volume = {31},
  issn = {0278-3649},
  doi = {10.1177/0278364912444766},
  language = {en},
  number = {7},
  journal = {The International Journal of Robotics Research},
  author = {Pavone, Marco and Smith, Stephen L and Frazzoli, Emilio and Rus, Daniela},
  month = jun,
  year = {2012},
  pages = {839-854}
}

@book{hensher2007handbook,
  title={Handbook of transport modelling},
  author={Hensher, David A and Button, Kenneth J},
  year={2007},
  publisher={Emerald Group Publishing Limited}
}

@inproceedings{drchal2015data,
  title={Data driven validation framework for multi-agent activity-based models},
  author={Drchal, Jan and {\v{C}}ertick{\'y}, Michal and Jakob, Michal},
  booktitle={International Workshop on Multi-Agent Systems and Agent-Based Simulation},
  pages={55--67},
  year={2015},
  organization={Springer}
}

@Article{jass2016,
  author={Jan Drchal and {\v{C}}ertick{\'y} and Michal Jakob},
  title={{VALFRAM: Validation Framework for Activity-Based Models}},
  journal={Journal of Artificial Societies and Social Simulation},
  year=2016,
  volume={19},
  number={3},
  pages={1-5},
  month={},
  doi={},
  url={https://ideas.repec.org/a/jas/jasssj/2015-85-3.html}
}

@article{alonso-moraOndemandHighcapacityRidesharing2017,
  title = {On-Demand High-Capacity Ride-Sharing via Dynamic Trip-Vehicle Assignment},
  author = {Alonso-Mora, Javier and Samaranayake, Samitha and Wallar, Alex and Frazzoli, Emilio and Rus, Daniela},
  date = {2017-01-17},
  journaltitle = {Proceedings of the National Academy of Sciences},
  shortjournal = {PNAS},
  volume = {114},
  number = {3},
  eprint = {28049820},
  eprinttype = {pmid},
  pages = {462--467},
  issn = {0027-8424, 1091-6490},
  doi = {10.1073/pnas.1611675114},
  url = {http://www.pnas.org/content/114/3/462},
  urldate = {2018-04-29},
  langid = {english}
}

@inproceedings{bischoffCitywideSharedTaxis2017,
  title = {City-Wide Shared Taxis: {{A}} Simulation Study in {{Berlin}}},
  shorttitle = {City-Wide Shared Taxis},
  booktitle = {2017 {{IEEE}} 20th {{International Conference}} on {{Intelligent Transportation Systems}} ({{ITSC}})},
  author = {Bischoff, Joschka and Maciejewski, Michal and Nagel, Kai},
  date = {2017-10},
  pages = {275--280},
  issn = {2153-0017},
  doi = {10.1109/ITSC.2017.8317926},
  eventtitle = {2017 {{IEEE}} 20th {{International Conference}} on {{Intelligent Transportation Systems}} ({{ITSC}})}
}

@article{bischoffSimulationCitywideReplacement2016,
  title = {Simulation of {{City-wide Replacement}} of {{Private Cars}} with {{Autonomous Taxis}} in {{Berlin}}},
  author = {Bischoff, Joschka and Maciejewski, Michal},
  date = {2016-01-01},
  journaltitle = {Procedia Computer Science},
  shortjournal = {Procedia Computer Science},
  series = {The 7th {{International Conference}} on {{Ambient Systems}}, {{Networks}} and {{Technologies}} ({{ANT}} 2016) / {{The}} 6th {{International Conference}} on {{Sustainable Energy Information Technology}} ({{SEIT-2016}}) / {{Affiliated Workshops}}},
  volume = {83},
  pages = {237--244},
  issn = {1877-0509},
  doi = {10.1016/j.procs.2016.04.121},
  url = {http://www.sciencedirect.com/science/article/pii/S1877050916301442},
  urldate = {2020-01-20},
  langid = {english}
}

@article{campbellEfficientInsertionHeuristics2004,
  title = {Efficient {{Insertion Heuristics}} for {{Vehicle Routing}} and {{Scheduling Problems}}},
  author = {Campbell, Ann and Savelsbergh, Martin},
  date = {2004-08-01},
  journaltitle = {Transportation Science},
  volume = {38},
  pages = {369--378},
  doi = {10.1287/trsc.1030.0046}
}

@inproceedings{capMultiObjectiveAnalysisRidesharing2018,
  title = {Multi-{{Objective Analysis}} of {{Ridesharing}} in {{Automated Mobility-on-Demand}}},
  booktitle = {Robotics: {{Science}} and {{Systems XIV}}},
  author = {Čáp, Michal and Alonso-Mora, Javier},
  date = {2018-06-26},
  publisher = {{Robotics: Science and Systems Foundation}},
  doi = {10.15607/RSS.2018.XIV.039},
  url = {http://www.roboticsproceedings.org/rss14/p39.pdf},
  urldate = {2019-01-18},
  eventtitle = {Robotics: {{Science}} and {{Systems}} 2018},
  isbn = {978-0-9923747-4-7},
  langid = {english}
}

@article{cordeauDialarideProblemModels2007,
  title = {The Dial-a-Ride Problem: Models and Algorithms},
  shorttitle = {The Dial-a-Ride Problem},
  author = {Cordeau, Jean-François and Laporte, Gilbert},
  date = {2007-09-01},
  journaltitle = {Annals of Operations Research},
  shortjournal = {Ann Oper Res},
  volume = {153},
  number = {1},
  pages = {29--46},
  issn = {1572-9338},
  doi = {10.1007/s10479-007-0170-8},
  url = {https://doi.org/10.1007/s10479-007-0170-8},
  urldate = {2019-01-12},
  langid = {english}
}

@article{drchalDatadrivenActivityScheduler2019,
  title = {Data-Driven Activity Scheduler for Agent-Based Mobility Models},
  author = {Drchal, Jan and Čertický, Michal and Jakob, Michal},
  date = {2019-01-01},
  journaltitle = {Transportation Research Part C: Emerging Technologies},
  shortjournal = {Transportation Research Part C: Emerging Technologies},
  volume = {98},
  pages = {370--390},
  issn = {0968-090X},
  doi = {10.1016/j.trc.2018.12.002},
  url = {https://www.sciencedirect.com/science/article/pii/S0968090X18306417},
  urldate = {2021-03-10},
  langid = {english}
}

@inproceedings{fiedlerImpactMobilityondemandTraffic2017,
  title = {Impact of Mobility-on-Demand on Traffic Congestion: {{Simulation-based}} Study},
  shorttitle = {Impact of Mobility-on-Demand on Traffic Congestion},
  booktitle = {2017 {{IEEE}} 20th {{International Conference}} on {{Intelligent Transportation Systems}} ({{ITSC}})},
  author = {Fiedler, D. and Čáp, M. and Čertický, M.},
  date = {2017-10},
  pages = {1--6},
  doi = {10.1109/ITSC.2017.8317830},
  eventtitle = {2017 {{IEEE}} 20th {{International Conference}} on {{Intelligent Transportation Systems}} ({{ITSC}})}
}

@inproceedings{fiedlerImpactRidesharingMobilityonDemand2018,
  title = {The {{Impact}} of {{Ridesharing}} in {{Mobility-on-Demand Systems}}: {{Simulation Case Study}} in {{Prague}}},
  shorttitle = {The {{Impact}} of {{Ridesharing}} in {{Mobility-on-Demand Systems}}},
  booktitle = {2018 21st {{International Conference}} on {{Intelligent Transportation Systems}} ({{ITSC}})},
  author = {Fiedler, David and Čertický, Michal and Alonso-Mora, Javier and Čáp, Michal},
  date = {2018-11},
  pages = {1173--1178},
  issn = {2153-0017, 2153-0009},
  doi = {10.1109/ITSC.2018.8569451},
  eventtitle = {2018 21st {{International Conference}} on {{Intelligent Transportation Systems}} ({{ITSC}})}
}

@article{fielbaumOndemandRidesharingOptimized2021,
  title = {On-Demand Ridesharing with Optimized Pick-up and Drop-off Walking Locations},
  author = {Fielbaum, Andres and Bai, Xiaoshan and Alonso-Mora, Javier},
  date = {2021-05-01},
  journaltitle = {Transportation Research Part C: Emerging Technologies},
  shortjournal = {Transportation Research Part C: Emerging Technologies},
  volume = {126},
  pages = {103061},
  issn = {0968-090X},
  doi = {10.1016/j.trc.2021.103061},
  url = {https://www.sciencedirect.com/science/article/pii/S0968090X21000887},
  urldate = {2022-06-04},
  langid = {english}
}

@article{fielbaumOptimizingVehicleRoute2021,
  title = {Optimizing a Vehicle’s Route in an on-Demand Ridesharing System in Which Users Might Walk},
  author = {Fielbaum, Andrés},
  date = {2021-03-18},
  journaltitle = {Journal of Intelligent Transportation Systems},
  volume = {0},
  number = {0},
  pages = {1--20},
  publisher = {{Taylor \& Francis}},
  issn = {1547-2450},
  doi = {10.1080/15472450.2021.1901225},
  url = {https://doi.org/10.1080/15472450.2021.1901225},
  urldate = {2022-06-04}
}

@article{hoSurveyDialarideProblems2018,
  title = {A Survey of Dial-a-Ride Problems: {{Literature}} Review and Recent Developments},
  shorttitle = {A Survey of Dial-a-Ride Problems},
  author = {Ho, Sin C. and Szeto, W. Y. and Kuo, Yong-Hong and Leung, Janny M. Y. and Petering, Matthew and Tou, Terence W. H.},
  date = {2018-05-01},
  journaltitle = {Transportation Research Part B: Methodological},
  shortjournal = {Transportation Research Part B: Methodological},
  volume = {111},
  pages = {395--421},
  issn = {0191-2615},
  doi = {10.1016/j.trb.2018.02.001},
  url = {http://www.sciencedirect.com/science/article/pii/S0191261517304484},
  urldate = {2019-01-14}
}

@article{jungDynamicSharedTaxiDispatch2015,
  title = {Dynamic {{Shared-Taxi Dispatch Algorithm}} with {{Hybrid Simulated Annealing}}},
  author = {Jung, Jaeyoung and Jayakrishnan, R and Young Park, Ji},
  date = {2015-06-22},
  journaltitle = {Computer-Aided Civil and Infrastructure Engineering},
  volume = {31},
  number = {4},
  pages = {275--291},
  doi = {10.1111/mice.12157}
}

@article{kalinaAgentsVehicleRouting2015,
  title = {Agents {{Toward Vehicle Routing Problem With Time Windows}}},
  author = {Kalina, Petr and Vokřínek, Jiří and Mařík, Vladimír},
  date = {2015-01-02},
  journaltitle = {Journal of Intelligent Transportation Systems},
  shortjournal = {Journal of Intelligent Transportation Systems},
  volume = {19},
  number = {1},
  pages = {3--17},
  publisher = {{Taylor \& Francis}},
  issn = {1547-2450},
  doi = {10.1080/15472450.2014.889953},
  url = {https://www.tandfonline.com/doi/10.1080/15472450.2014.889953},
  urldate = {2020-06-16}
}

@article{liRestrictedPathbasedRidesharing2019,
  title = {A Restricted Path-Based Ridesharing User Equilibrium},
  author = {Li, Meng and Di, Xuan and Liu, Henry X. and Huang, Hai-Jun},
  date = {2019-09-13},
  journaltitle = {Journal of Intelligent Transportation Systems},
  volume = {0},
  number = {0},
  pages = {1--21},
  publisher = {{Taylor \& Francis}},
  issn = {1547-2450},
  doi = {10.1080/15472450.2019.1658525},
  url = {https://doi.org/10.1080/15472450.2019.1658525},
  urldate = {2020-05-21}
}

@article{maciejewskiCongestionEffectsAutonomous2018,
  title = {Congestion Effects of Autonomous Taxi Fleets},
  author = {Maciejewski, Michal and Bischoff, Joschka},
  date = {2018-12-05},
  journaltitle = {Transport},
  shortjournal = {1},
  volume = {33},
  number = {4},
  pages = {971--980},
  issn = {1648-3480},
  doi = {10.3846/16484142.2017.1347827},
  url = {https://journals.vgtu.lt/index.php/Transport/article/view/212},
  urldate = {2019-03-28},
  langid = {english}
}

@article{maDynamicRidesharingDispatch2019,
  title = {A Dynamic Ridesharing Dispatch and Idle Vehicle Repositioning Strategy with Integrated Transit Transfers},
  author = {Ma, Tai-Yu and Rasulkhani, Saeid and Chow, Joseph Y. J. and Klein, Sylvain},
  date = {2019-08-01},
  journaltitle = {Transportation Research Part E: Logistics and Transportation Review},
  shortjournal = {Transportation Research Part E: Logistics and Transportation Review},
  volume = {128},
  pages = {417--442},
  issn = {1366-5545},
  doi = {10.1016/j.tre.2019.07.002},
  url = {http://www.sciencedirect.com/science/article/pii/S1366554518314790},
  urldate = {2020-02-06},
  langid = {english}
}

@article{maRidesharingUserEquilibrium2020,
  title = {Ridesharing User Equilibrium Problem under {{OD-based}} Surge Pricing Strategy},
  author = {Ma, Jie and Xu, Min and Meng, Qiang and Cheng, Lin},
  date = {2020-04-01},
  journaltitle = {Transportation Research Part B: Methodological},
  shortjournal = {Transportation Research Part B: Methodological},
  volume = {134},
  pages = {1--24},
  issn = {0191-2615},
  doi = {10.1016/j.trb.2020.02.001},
  url = {https://www.sciencedirect.com/science/article/pii/S0191261519303832},
  urldate = {2021-11-27},
  langid = {english}
}

@article{masmoudiThreeEffectiveMetaheuristics2016,
  title = {Three Effective Metaheuristics to Solve the Multi-Depot Multi-Trip Heterogeneous Dial-a-Ride Problem},
  author = {Masmoudi, Mohamed Amine and Hosny, Manar and Braekers, Kris and Dammak, Abdelaziz},
  date = {2016-12-01},
  journaltitle = {Transportation Research Part E: Logistics and Transportation Review},
  shortjournal = {Transportation Research Part E: Logistics and Transportation Review},
  volume = {96},
  pages = {60--80},
  issn = {1366-5545},
  doi = {10.1016/j.tre.2016.10.002},
  url = {http://www.sciencedirect.com/science/article/pii/S1366554516304070},
  urldate = {2020-02-07},
  langid = {english}
}

@article{masoudRealtimeAlgorithmSolve2017,
  title = {A Real-Time Algorithm to Solve the Peer-to-Peer Ride-Matching Problem in a Flexible Ridesharing System},
  author = {Masoud, Neda and Jayakrishnan, R.},
  date = {2017-12-01},
  journaltitle = {Transportation Research Part B: Methodological},
  shortjournal = {Transportation Research Part B: Methodological},
  volume = {106},
  pages = {218--236},
  issn = {0191-2615},
  doi = {10.1016/j.trb.2017.10.006},
  url = {http://www.sciencedirect.com/science/article/pii/S0191261517301169},
  urldate = {2020-02-06},
  langid = {english}
}

@inproceedings{millerPredictivePositioningQuality2017,
  title = {Predictive Positioning and Quality of Service Ridesharing for Campus Mobility on Demand Systems},
  booktitle = {2017 {{IEEE International Conference}} on {{Robotics}} and {{Automation}} ({{ICRA}})},
  author = {Miller, J. and How, J. P.},
  date = {2017-05},
  pages = {1402--1408},
  doi = {10.1109/ICRA.2017.7989167},
  eventtitle = {2017 {{IEEE International Conference}} on {{Robotics}} and {{Automation}} ({{ICRA}})}
}

@article{muelasDistributedVNSAlgorithm2015,
  title = {A Distributed {{VNS}} Algorithm for Optimizing Dial-a-Ride Problems in Large-Scale Scenarios},
  author = {Muelas, Santiago and LaTorre, Antonio and Peña, José-María},
  date = {2015-05-01},
  journaltitle = {Transportation Research Part C: Emerging Technologies},
  shortjournal = {Transportation Research Part C: Emerging Technologies},
  volume = {54},
  pages = {110--130},
  issn = {0968-090X},
  doi = {10.1016/j.trc.2015.02.024},
  url = {http://www.sciencedirect.com/science/article/pii/S0968090X15000790},
  urldate = {2020-02-07},
  langid = {english}
}

@article{muelasVariableNeighborhoodSearch2013,
  title = {A Variable Neighborhood Search Algorithm for the Optimization of a Dial-a-Ride Problem in a Large City},
  author = {Muelas, Santiago and LaTorre, Antonio and Peña, José-María},
  date = {2013-10-15},
  journaltitle = {Expert Systems with Applications},
  shortjournal = {Expert Systems with Applications},
  volume = {40},
  number = {14},
  pages = {5516--5531},
  issn = {0957-4174},
  doi = {10.1016/j.eswa.2013.04.015},
  url = {http://www.sciencedirect.com/science/article/pii/S0957417413002522},
  urldate = {2020-02-04},
  langid = {english}
}

@report{nyctaxilimousinecommission2016TLCFactbook2016,
  title = {2016 {{TLC Factbook}}},
  author = {{NYC Taxi \& Limousine Commission}},
  date = {2016},
  pages = {15},
  institution = {{NYC Taxi \& Limousine Commission}}
}

@report{nyctaxilimousinecommission2018Factbook2018,
  title = {2018 {{Factbook}}},
  author = {{NYC Taxi \& Limousine Commission}},
  date = {2018},
  pages = {15},
  institution = {{NYC Taxi \& Limousine Commission}}
}

@article{pavoneRoboticLoadBalancing2012,
  title = {Robotic Load Balancing for Mobility-on-Demand Systems},
  author = {Pavone, Marco and Smith, Stephen L and Frazzoli, Emilio and Rus, Daniela},
  date = {2012-06-01},
  journaltitle = {The International Journal of Robotics Research},
  shortjournal = {The International Journal of Robotics Research},
  volume = {31},
  number = {7},
  pages = {839--854},
  issn = {0278-3649},
  doi = {10.1177/0278364912444766},
  url = {https://doi.org/10.1177/0278364912444766},
  urldate = {2019-01-26},
  langid = {english}
}

@article{pfeifferALNSAlgorithmStatic2022,
  title = {An {{ALNS}} Algorithm for the Static Dial-a-Ride Problem with Ride and Waiting Time Minimization},
  author = {Pfeiffer, Christian and Schulz, Arne},
  date = {2022-03-01},
  journaltitle = {OR Spectrum},
  shortjournal = {OR Spectrum},
  volume = {44},
  number = {1},
  pages = {87--119},
  issn = {1436-6304},
  doi = {10.1007/s00291-021-00656-7},
  url = {https://doi.org/10.1007/s00291-021-00656-7},
  urldate = {2022-06-04},
  langid = {english}
}

@article{santiQuantifyingBenefitsVehicle2014,
  title = {Quantifying the Benefits of Vehicle Pooling with Shareability Networks},
  author = {Santi, Paolo and Resta, Giovanni and Szell, Michael and Sobolevsky, Stanislav and Strogatz, Steven H. and Ratti, Carlo},
  date = {2014-09-16},
  journaltitle = {Proceedings of the National Academy of Sciences},
  volume = {111},
  number = {37},
  pages = {13290--13294},
  issn = {0027-8424, 1091-6490},
  doi = {10.1073/pnas.1403657111},
  url = {http://www.pnas.org/lookup/doi/10.1073/pnas.1403657111},
  urldate = {2019-01-26},
  langid = {english}
}

@inproceedings{santosDynamicTaxiRidesharing2013,
  title = {Dynamic {{Taxi}} and {{Ridesharing}}: {{A Framework}} and {{Heuristics}} for the {{Optimization Problem}}},
  shorttitle = {Dynamic {{Taxi}} and {{Ridesharing}}},
  booktitle = {Proceedings of the {{Twenty-Third International Joint Conference}} on {{Artificial Intelligence}}},
  author = {Santos, Douglas O. and Xavier, Eduardo C.},
  date = {2013},
  series = {{{IJCAI}} '13},
  pages = {2885--2891},
  publisher = {{AAAI Press}},
  url = {http://dl.acm.org/citation.cfm?id=2540128.2540544},
  urldate = {2019-03-28},
  isbn = {978-1-57735-633-2},
  venue = {Beijing, China}
}

@article{santosTaxiRideSharing2015,
  title = {Taxi and {{Ride Sharing}}: {{A Dynamic Dial-a-Ride Problem}} with {{Money}} as an {{Incentive}}},
  shorttitle = {Taxi and {{Ride Sharing}}},
  author = {Santos, Douglas O. and Xavier, Eduardo C.},
  date = {2015-11-01},
  journaltitle = {Expert Systems with Applications},
  shortjournal = {Expert Systems with Applications},
  volume = {42},
  number = {19},
  pages = {6728--6737},
  issn = {0957-4174},
  doi = {10.1016/j.eswa.2015.04.060},
  url = {https://www.sciencedirect.com/science/article/pii/S0957417415003024},
  urldate = {2021-03-02},
  langid = {english}
}

@book{schrijverTheoryLinearInteger1986,
  title = {Theory of Linear and Integer Programming},
  author = {Schrijver, Alexander},
  date = {1986},
  publisher = {{John Wiley \& Sons, Inc.}},
  location = {{USA}},
  isbn = {978-0-471-90854-8},
  pagetotal = {471}
}

@incollection{spieserSystematicApproachDesign2014a,
  title = {Toward a {{Systematic Approach}} to the {{Design}} and {{Evaluation}} of {{Automated Mobility-on-Demand Systems}}: {{A Case Study}} in {{Singapore}}},
  shorttitle = {Toward a {{Systematic Approach}} to the {{Design}} and {{Evaluation}} of {{Automated Mobility-on-Demand Systems}}},
  booktitle = {Road {{Vehicle Automation}}},
  author = {Spieser, Kevin and Treleaven, Kyle and Zhang, Rick and Frazzoli, Emilio and Morton, Daniel and Pavone, Marco},
  editor = {Meyer, Gereon and Beiker, Sven},
  date = {2014},
  series = {Lecture {{Notes}} in {{Mobility}}},
  pages = {229--245},
  publisher = {{Springer International Publishing}},
  location = {{Cham}},
  doi = {10.1007/978-3-319-05990-7_20},
  url = {https://doi.org/10.1007/978-3-319-05990-7_20},
  urldate = {2020-05-11},
  isbn = {978-3-319-05990-7},
  langid = {english}
}

@incollection{spieserSystematicApproachDesign2014b,
  title = {Toward a {{Systematic Approach}} to the {{Design}} and {{Evaluation}} of {{Automated Mobility-on-Demand Systems}}: {{A Case Study}} in {{Singapore}}},
  shorttitle = {Toward a {{Systematic Approach}} to the {{Design}} and {{Evaluation}} of {{Automated Mobility-on-Demand Systems}}},
  booktitle = {Road {{Vehicle Automation}}},
  author = {Spieser, Kevin and Treleaven, Kyle and Zhang, Rick and Frazzoli, Emilio and Morton, Daniel and Pavone, Marco},
  editor = {Meyer, Gereon and Beiker, Sven},
  date = {2014},
  series = {Lecture {{Notes}} in {{Mobility}}},
  pages = {229--245},
  publisher = {{Springer International Publishing}},
  location = {{Cham}},
  doi = {10.1007/978-3-319-05990-7_20},
  url = {https://doi.org/10.1007/978-3-319-05990-7_20},
  urldate = {2021-12-08},
  isbn = {978-3-319-05990-7},
  langid = {english}
}

@article{susanshaheenSimilaritiesDifferencesMobility2020,
  title = {Similarities and {{Differences}} of {{Mobility}} on {{Demand}} ({{MOD}}) and {{Mobility}} as a {{Service}} ({{MaaS}}) | {{Transportation Sustainability Research Center}}},
  author = {{Susan Shaheen} and {Adam Cohen}},
  date = {2020-06-01},
  journaltitle = {ite journal},
  volume = {90},
  number = {6},
  url = {https://tsrc.berkeley.edu/publications/similarities-and-differences-mobility-demand-mod-and-mobility-service-maas},
  urldate = {2021-09-09},
  langid = {english}
}

@inproceedings{tadakiCriticalDensityExperimental2015a,
  title = {Critical {{Density}} of {{Experimental Traffic Jam}}},
  booktitle = {Traffic and {{Granular Flow}} '13},
  author = {Tadaki, Shin-ichi and Kikuchi, Macoto and Fukui, Minoru and Nakayama, Akihiro and Nishinari, Katsuhiro and Shibata, Akihiro and Sugiyama, Yuki and Yosida, Taturu and Yukawa, Satoshi},
  editor = {Chraibi, Mohcine and Boltes, Maik and Schadschneider, Andreas and Seyfried, Armin},
  date = {2015},
  pages = {505--511},
  publisher = {{Springer International Publishing}},
  location = {{Cham}},
  doi = {10.1007/978-3-319-10629-8_56},
  isbn = {978-3-319-10629-8},
  langid = {english}
}

@article{tamannaeiCarpoolingProblemNew2019,
  title = {Carpooling Problem: {{A}} New Mathematical Model, Branch-and-Bound, and Heuristic Beam Search Algorithm},
  shorttitle = {Carpooling Problem},
  author = {Tamannaei, Mohammad and Irandoost, Iman},
  date = {2019-05-04},
  journaltitle = {Journal of Intelligent Transportation Systems},
  volume = {23},
  number = {3},
  pages = {203--215},
  publisher = {{Taylor \& Francis}},
  issn = {1547-2450},
  doi = {10.1080/15472450.2018.1484739},
  url = {https://doi.org/10.1080/15472450.2018.1484739},
  urldate = {2020-05-21}
}

@book{tothVehicleRoutingProblems2014,
  title = {Vehicle {{Routing}}: {{Problems}}, {{Methods}}, and {{Applications}}, {{Second Edition}}},
  shorttitle = {Vehicle {{Routing}}},
  author = {Toth, Paolo and Vigo, Daniele},
  date = {2014-12-05},
  eprint = {AoTTBQAAQBAJ},
  eprinttype = {googlebooks},
  publisher = {{SIAM}},
  isbn = {978-1-61197-359-4},
  langid = {english},
  pagetotal = {467}
}

@article{vanengelenEnhancingFlexibleTransport2018,
  title = {Enhancing Flexible Transport Services with Demand-Anticipatory Insertion Heuristics},
  author = {van Engelen, Matti and Cats, Oded and Post, Henk and Aardal, Karen},
  options = {useprefix=true},
  date = {2018-02-01},
  journaltitle = {Transportation Research Part E: Logistics and Transportation Review},
  shortjournal = {Transportation Research Part E: Logistics and Transportation Review},
  volume = {110},
  pages = {110--121},
  issn = {1366-5545},
  doi = {10.1016/j.tre.2017.12.015},
  url = {http://www.sciencedirect.com/science/article/pii/S1366554517307810},
  urldate = {2020-02-06},
  langid = {english}
}

@article{venkatramanCongestionawareTabuSearch2019,
  title = {A Congestion-Aware {{Tabu}} Search Heuristic to Solve the Shared Autonomous Vehicle Routing Problem},
  author = {Venkatraman, Prashanth and Levin, Michael W.},
  date = {2019-09-27},
  journaltitle = {Journal of Intelligent Transportation Systems},
  volume = {0},
  number = {0},
  pages = {1--13},
  publisher = {{Taylor \& Francis}},
  issn = {1547-2450},
  doi = {10.1080/15472450.2019.1665521},
  url = {https://doi.org/10.1080/15472450.2019.1665521},
  urldate = {2020-06-07}
}

@inproceedings{wallarOptimizingVehicleDistributions2019a,
  title = {Optimizing {{Vehicle Distributions}} and {{Fleet Sizes}} for {{Shared Mobility-on-Demand}}},
  booktitle = {2019 {{International Conference}} on {{Robotics}} and {{Automation}} ({{ICRA}})},
  author = {Wallar, A. and Alonso-Mora, J. and Rus, D.},
  date = {2019-05},
  pages = {3853--3859},
  issn = {2577-087X},
  doi = {10.1109/ICRA.2019.8793685},
  eventtitle = {2019 {{International Conference}} on {{Robotics}} and {{Automation}} ({{ICRA}})}
}

@article{yanStochasticRidesharingUser2019,
  title = {Stochastic {{Ridesharing User Equilibrium}} in {{Transport Networks}}},
  author = {Yan, Chen-Yang and Hu, Mao-Bin and Jiang, Rui and Long, Jiancheng and Chen, Jin-Yong and Liu, Hao-Xiang},
  date = {2019-12-01},
  journaltitle = {Networks and Spatial Economics},
  shortjournal = {Netw Spat Econ},
  volume = {19},
  number = {4},
  pages = {1007--1030},
  issn = {1572-9427},
  doi = {10.1007/s11067-019-9442-5},
  url = {https://doi.org/10.1007/s11067-019-9442-5},
  urldate = {2021-11-27},
  langid = {english}
}

@article{zhanModifiedArtificialBee2021,
  title = {A Modified Artificial Bee Colony Algorithm for the Dynamic Ride-Hailing Sharing Problem},
  author = {Zhan, Xingbin and Szeto, W. Y. and Shui, C. S. and Chen, Xiqun},
  date = {2021-06-01},
  journaltitle = {Transportation Research Part E: Logistics and Transportation Review},
  shortjournal = {Transportation Research Part E: Logistics and Transportation Review},
  volume = {150},
  pages = {102124},
  issn = {1366-5545},
  doi = {10.1016/j.tre.2020.102124},
  url = {https://www.sciencedirect.com/science/article/pii/S1366554520307729},
  urldate = {2022-06-04},
  langid = {english}
}

@article{zhanSimulationOptimizationFramework2022,
  title = {A Simulation–Optimization Framework for a Dynamic Electric Ride-Hailing Sharing Problem with a Novel Charging Strategy},
  author = {Zhan, Xingbin and Szeto, W. Y. and (Michael) Chen, Xiqun},
  date = {2022-03-01},
  journaltitle = {Transportation Research Part E: Logistics and Transportation Review},
  shortjournal = {Transportation Research Part E: Logistics and Transportation Review},
  volume = {159},
  pages = {102615},
  issn = {1366-5545},
  doi = {10.1016/j.tre.2022.102615},
  url = {https://www.sciencedirect.com/science/article/pii/S1366554522000138},
  urldate = {2022-06-04},
  langid = {english}
}

\begin{figure}
    \centering
    \includegraphics[width=0.5\columnwidth]{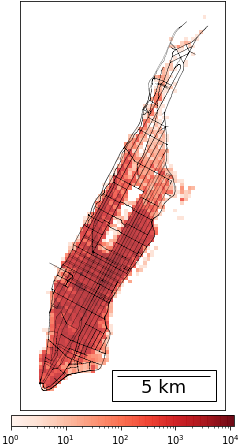}
    \caption{\rev{Manhattan taxi trip requests on Friday, May 10, 2013, between 19:00 and 20:00. The start positions of all vehicle trips are discretized to squares of \num{200} square meters. Darker color translates to higher demand, and the color bar has a logarithmic scale.}}
    \label{fig:demand_manhattan}
\end{figure}

\appendix
\rev{

\section{Manhattan Experiment}
\label{sec:manhattan}
In order to demonstrate general applicability of our method and to allow for easier comparison with previous work, we repeat our experiment in Manhattan using a publicly available dataset of taxi trips as transportation demand.
Specifically, we use the same demand and road network as used by~\textcite{alonso-moraOndemandHighcapacityRidesharing2017}.
Identically to our Prague experiment, we simulated the system for one hour with a 30-minute warm-up period. 
While~\textcite{alonso-moraOndemandHighcapacityRidesharing2017} run the simulation for one week worth of data, here, for simplicity, we selected the day and hour with the largest number of requests, which was Friday, May 10, 2013, between 19:00 and 20:00.
\srev{There are \num{137202} travel request in the selected period, }Figure~\ref{fig:demand_manhattan} shows the spatial structure of the demand.
}

\begin{figure}
    \centering
    \includegraphics[width=0.5\columnwidth]{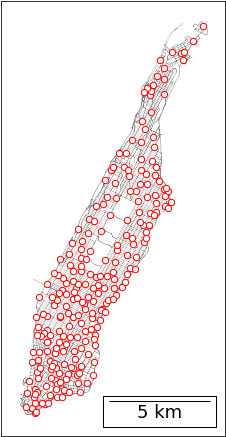}
    \caption{\rev{MoD stations on Manhattan. \srev{Each red circle represents a single MoD system station.}}}
    \label{fig:stations_manhattan}
\end{figure}

\rev{
Like~\textcite{alonso-moraOndemandHighcapacityRidesharing2017}, we use travel speeds along individual road segments derived from historical data, but instead of computing the speeds from the travel demand, we use the speeds from the Uber Movement\footnote{\url{https://movement.uber.com/}} open data project.
Other than that, we followed the methodology described in the main part of this article. 
That is, we assume a station-based model and perform rebalancing, fleet sizing, and passenger-vehicle matching as described in Section~\ref{sec:methodology}.
Because the historical speeds from Uber Movement dataset are on average approximately half of the posted speed and there are a lot of one-way streets on Manhattan, we need \num{236} stations to provide the required quality of service, even though Manhattan is about five times smaller than Prague. 
Figure~\ref{fig:stations_manhattan} shows the locations of stations.
}

\rev{
On Manhattan, we evaluated five of the six scenarios tested in the Prague case study. 
We do not evaluate the present state scenario, as the Manhattan dataset represents taxi trips, and therefore, the scenario with MoD system without ridesharing is, in fact, also the ``present state'' scenario.
}
\srev{
In Table~\ref{fig:comparison_table-manhattan}, we can see the results of the same set of experiments as we performed in the Prague case study.
Our optimal implementation of the VGA method was able to compute the optimal assignments while the average computational time for a \SI{30}{seconds} batch was less than \SI{7}{seconds}.
This is in contrast to results reported in~\textcite{alonso-moraOndemandHighcapacityRidesharing2017} that were not computed to optimality and required more than \SI{21}{seconds} to compute the most similar configuration ($ q_{\mathrm{max}} =$ \SI{5}{minutes}, vehicle capacity of four passengers, \SI{3000}{vehicles}). 
This may be because the algorithm by~\textcite{alonso-moraOndemandHighcapacityRidesharing2017} was developed and optimized to allow evaluation of scenarios with even larger delays of 7 minutes and with vehicle capacities of up to 10 passengers; such configurations result in an exponentially larger number of potential passenger-vehicle assignments and consequently, cannot be computed to optimality even with our performance-optimized VGA method.
}

\begin{table}
\scriptsize
\centering{}%
\setlength{\tabcolsep}{0.3em}
{\renewcommand{\arraystretch}{1.2}%
\begin{tabular}{|+l|-r|-r|-r|-r|-r|}
\hline
 & \thead{No Ridesh.} & \thead{IH} & \thead{VGA} & \thead{VGA lim} & \thead{VGA PNAS*}
\tabularnewline
\hline
\hline
Optimal & - & no & yes & no & no
\tabularnewline
\hline
\rowstyle{\bfseries}
Total veh. dist. (km) & \num{868899} & \num{362387} & \num{334195} & \num{334737}& \num{377563} (\num{344057})
\tabularnewline
\hline
Avg. delay (s) & \num{109} & \num{117} & \num{109} & \num{109} & \num{83} (\num{110})
\tabularnewline
\hline
Avg. density (veh/km) & \num{0.0183} & \num{0.0085} & \num{0.008} & \num{0.008} & \num{0.0089} (\num{0.0082})
\tabularnewline
\hline
Congested seg. & \num{220} & \num{9} & \num{9} & \num{8} & \num{14} (\num{10})
\tabularnewline
\hline
Heavily loaded seg.  & \num{692} & \num{129} & \num{117} & \num{112} & \num{152} (\num{115})
\tabularnewline
\hline
Used Vehicles  & \num{46186} & \num{20272} & \num{19714} & \num{19712} & \num{22545} (\num{20293})
\tabularnewline
\hline
Avg. comp. time (ms)  & \num{918} & \num{57} & \num{6646} & \num{6650} & \num{24717} (\num{7170})
\tabularnewline
\hline
\end{tabular}}

\caption{\label{fig:comparison_table-manhattan} \rev{
Main results from the Manhattan scenarios during the peak (19:00-20:00) with a maximum passenger delay of \SI{4}{minutes}. 
Congested segments are segments on which traffic density is above critical density, and heavily loaded segments are segments with a density above \SI{50}{\percent} of the critical density.
For the VGA PNAS method, we also tested a version that does not limit the vehicles considered for each request to 30 nearest vehicles (in parentheses).
}}
\end{table}

\rev{
Because the Manhattan experiment is less complex compared to the Prague experiment, we can observe a similar effect as in the Prague off-peak experiment: the VGA limited method computes only slightly worse solutions than the optimal method, and also the computational times are similar. 
This is because the time limits of the VGA limited method were not reached in the majority of iterations.

Our re-implementation of the PNAS method gives a rather surprising result: the performance metrics are worse than the IH while using more computational time than the optimal method.
We investigated this surprising result and found out that the cause is one of the heuristics that limits the number of vehicles considered for assignment to a particular request to 30 nearest vehicles.
This heuristic can limit the exploration so much that the solution can be worse than the IH solution. 
Moreover, for less complex scenarios, the time needed to compute the 30 nearest vehicles can dominate the total computational time, as it happened in our case, probably because this heuristic was not optimized.
This observation suggests that in order to achieve acceptable performance, one may need to vary the parameters of heuristics based on the complexity of the problem instance at hand.
We also performed the experiment using the VGA PNAS method with this heuristic turned off. 
For results, see numbers in parentheses in the last column of result tables (\ref{fig:comparison_table-manhattan}, \ref{fig:comparison_table-manhattan_2}).
}

\srev{
Table~\ref{fig:comparison_table-manhattan_2} shows another set of results of experiments with $ q_{\mathrm{max}} = $ \SI{7}{minutes} and the capacity of 10 persons per vehicle, which corresponds to the most complex configuration in~\textcite{alonso-moraOndemandHighcapacityRidesharing2017}.
In this experiment set, we only evaluated the three sub-optimal ridesharing methods to see how they behave under such parametrization.
Interestingly, for this scenario, the IH achieves the best performance: The IH finds plans with a total traveled distance that is \SI{12}{\percent} smaller than plans found by both sub-optimal versions of the VGA method using only a fraction of computational resources.}
\trev{This experiment demonstrates the limit of applicability of the VGA method for routing in large-scale MoD systems. The relaxed time windows and increased vehicle capacity increase the number of feasible groups and the maximum group size to a level that cannot be solved by the ILP solver and the single-vehicle solver, respectively, in practical time. 
Consequently, the VGA algorithm is unable to return an optimal solution to such instances. 
}

\begin{table}
\scriptsize
\centering{}%
\setlength{\tabcolsep}{0.3em}
{\renewcommand{\arraystretch}{1.2}%
\begin{tabular}{|+l|-r|-r|-r|-r|-r|}
\hline
 & \thead{IH} & \thead{VGA lim} & \thead{VGA PNAS*}
\tabularnewline
\hline
\hline
Optimal & no & no & no
\tabularnewline
\hline
\rowstyle{\bfseries}
Total veh. dist. (km) & \num{233859} & \num{275028} & \num{267471}
\tabularnewline
\hline
Avg. delay (s) & \num{227} & \num{224} & \num{217}
\tabularnewline
\hline
Avg. density (veh/km) & \num{0.006} & \num{0.0067} & \num{0.0066}
\tabularnewline
\hline
Congested seg. & \num{0} & \num{1} & \num{0}
\tabularnewline
\hline
Heavily loaded seg.  & \num{24} & \num{53} & \num{48}
\tabularnewline
\hline
Used Vehicles  & \num{13319} & \num{16517} & \num{16025}
\tabularnewline
\hline
Avg. comp. time (ms)  & \num{25} & \num{57554} & \num{139811}
\tabularnewline
\hline
\end{tabular}}

\caption{\label{fig:comparison_table-manhattan_2}
\srev{Main results from the Manhattan scenarios during the peak (19:00-20:00) with a maximum passenger delay of \SI{7}{minutes} and vehicle capacity of 10 persons per vehicle. 
Congested segments are segments on which traffic density is above critical density, and heavily loaded segments are segments with a density above \SI{50}{\percent} of the critical density.
For the VGA PNAS method, we used the version that does \emph{not} limit the vehicles considered for each request to 30 nearest vehicles.
}}
\end{table}

\begin{figure}
\centering{}\includegraphics[width=0.9\columnwidth]{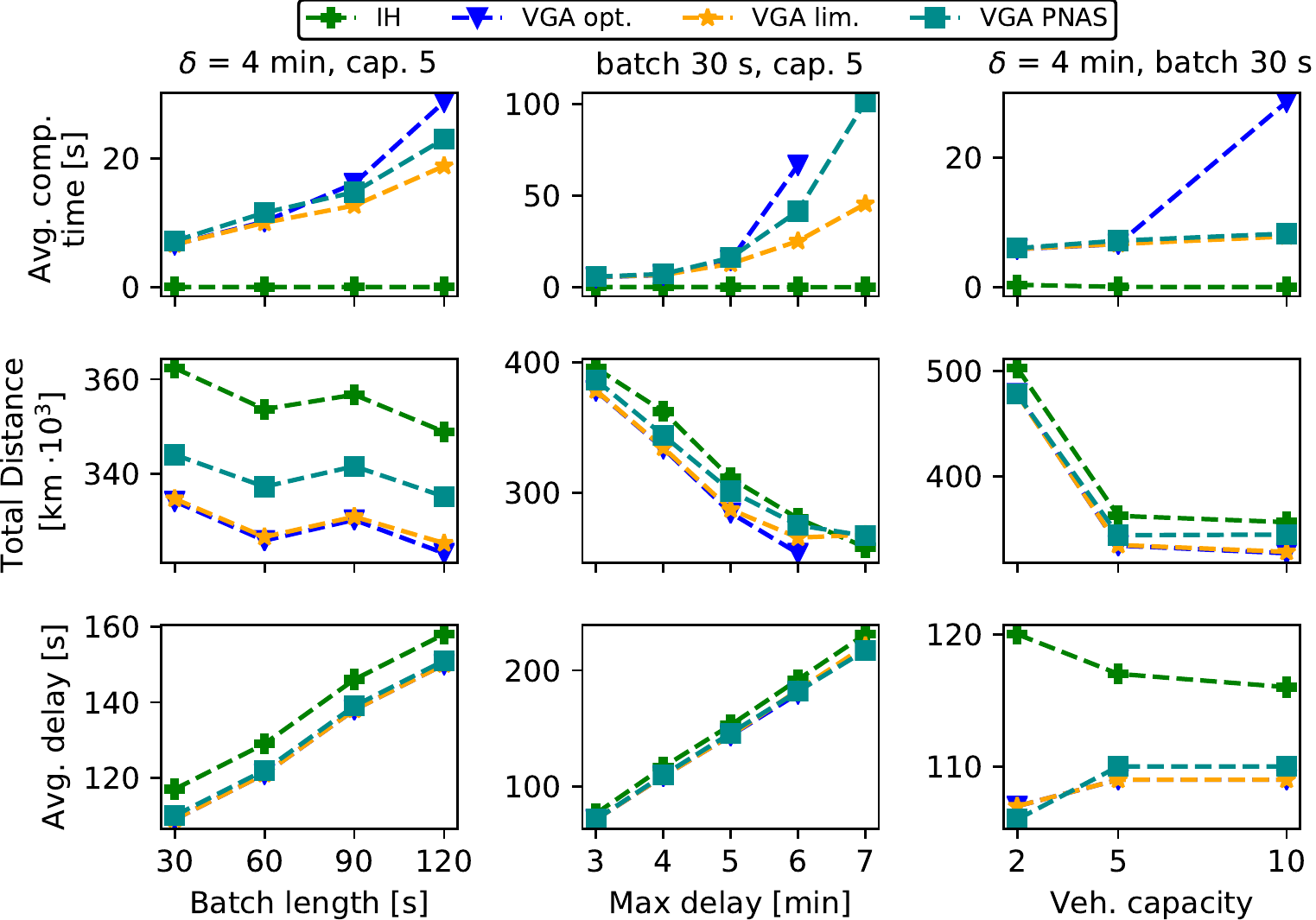}
\caption{\label{fig:sensitivity_analysis-manhattan} \srev{Sensitivity analysis: Manhattan. 
Each column represents one experiment set, and inside each column, each value on the x-axis represents one experiment. 
Each row displays a single measured quantity. 
The optimal VGA method is only computed for the maximum delay of up to 6 minutes.
For the maximum delay of 7 minutes, the optimal ridesharing assignment cannot be computed within 24 hours run time limit.}}
\end{figure}

\srev{Finally, we performed a sensitivity analysis for the Manhattan case study: the results are reported in Figure~\ref{fig:sensitivity_analysis-manhattan}.
It tells a similar story as the sensitivity analysis for the Prague case study (Section~\ref{sec:results}). 
Some of the previously discussed phenomena are even more apparent in the Manhattan sensitivity analysis. 
We can see that the computational requirements grow with the maximum delay not only for the optimal VGA method but also for the resource-constrained VGA methods.
Also, we can clearly see that the efficiency (total distance driven) gap between IH and the constrained VGA methods is shrinking for larger maximum delays.
For the maximum delay of 7 minutes, the IH method starts to outperform both resource constrained VGA methods.
}

\rev{
\section{VGA Optimizations}
\label{sec:vga_optimizatios}
The main objective of this article is to quantify the performance gap between optimal ridesharing assignments and assignments computed using heuristic solutions.
However, a naive implementation of the VGA algorithm would require prohibitively long computation time and an extreme amount of memory to compute an optimal solution.
In order to arrive at the optimal solutions in a manageable time, we had to implement several performance optimizations.
Without these optimizations, the VGA method would need several hours to compute an optimal ridesharing assignment for each 30-second-long batch.

To reduce the number of vehicles considered in request-vehicle matching, we leverage the specific properties of the station-based MoD system and modify the VGA method accordingly.
We reduce the number of vehicles for which the groups are generated as follows:
First, we observe that we need at most as many vehicles as the number of waiting requests since, in the worst case, each request can be transported in a dedicated vehicle from the nearest station.
Second, we exploit symmetries in the solution space. 
We observe that vehicles parked in a station can be arbitrarily relabeled without any effect on the solution quality.
Therefore, instead of computing feasible groups for each vehicle parked in a station, we generate feasible groups for only one vehicle from that station, representing any other vehicle currently parked in the station. 
Consequently, in the assignment ILP, we can relax Constraint~\ref{vga-constr:1} corresponding to this representative vehicle to allow assigning as many vehicle plans as there are vehicles parked in the station:
\begin{equation}
    \sum_{g = 1}^{|\FG|} \xi_v^g \leq |V_s| \quad \forall{s \in S}. 
\end{equation}
In this modified version of Constraint~\ref{vga-constr:1}, $ V_s $ is the set of vehicles parked in station $ s $, and $ S $ is the set of all stations.

We paid special attention to an efficient implementation of the function $ \feas $ that is used to determine the feasibility of a newly formed group and to compute the optimal route for the group.
This function solves a single-vehicle DARP by searching through all feasible permutations of travel schedules.
Most of the time during group generation is spent inside this function. 
We achieved significant performance gains by implementing the algorithm in a way that constructs permutations "in place" and avoids memory allocation during the search process.
Also, we implemented a look-ahead procedure that is triggered each time a pickup or drop-off location is added to extend a partial plan.
In this look-ahead, we verify that the maximum allowed time for each pickup and drop-off location that are still waiting to be added to the plan is higher than the time of the most recently added order to the partial plan; if the above does not hold, we can safely discard the partial plan as infeasible.


Finally, we parallelized the group generation process so that feasible groups are computed in a separate thread for each vehicle.
An even better approach would be to parallelize the function $ \feas $ because that way, we can distribute the work among threads even if there is only a small number of vehicles with a larger number of vehicle groups.
However, we leave this optimization for future work.

}

\end{document}